\documentclass[11pt]{article}
\usepackage{amsfonts}
\usepackage{latexsym}
\usepackage{cite}
\usepackage{amsmath,amsfonts,latexsym,amssymb}
\usepackage[mathscr]{eucal}
\usepackage{cases,color}
\usepackage{amsthm}
\usepackage[bf,small]{caption2}
\usepackage{float}
\usepackage{graphicx}
\usepackage{amsmath}
\usepackage{amssymb}
\usepackage[all]{xy}
\usepackage{hyperref}

\newtheorem{theorem}{theorem}[section]

\newtheorem{cla}[theorem]{Claim}

\newtheorem{lem}[theorem]{Lemma}

\newtheorem{prob}[theorem]{Problem}
\newtheorem{prop}[theorem]{Proposition}

\newtheorem{rmk}[theorem]{Remark}
\newtheorem{thm}[theorem]{Theorem}

\begin{document}

\title{\vspace{-2cm}\textbf{Torsion in Kauffman bracket skein module of a $4$-strand Montesinos knot exterior}}
\author{\Large Haimiao Chen}
\date{}
\maketitle

\begin{abstract}
  For an oriented $3$-manifold $M$, let $\mathcal{S}(M)$ denote its Kauffman bracket skein module over $\mathbb{Z}[q^{\pm\frac{1}{2}}]$. We show that $\mathcal{S}(M)$ admits torsion when $M$ is the exterior of the Montesinos knot $K(a_1/b_1,a_2/b_2,a_3/b_4,a_4/b_4)$ with each $b_i\ge 3$. This provides a negative answer to Problem 1.92 (G)-(i) in the Kirby's list, which asks whether $\mathcal{S}(M)$ is free when $M$ is irreducible and has no incompressible non-boundary parallel torus.

\end{abstract}

\section{Introduction}

Let $R$ be a commutative ring with a distinguished invertible element $q^{\frac{1}{2}}$.
Given an oriented 3-manifold $M$, its {\it Kauffman bracket skein module}, denoted by $\mathcal{S}(M;R)$, is defined as the quotient of the free $R$-module generated by isotopy classes of (possibly empty) framed links embedded in $M$ by the submodule generated by the {\it skein relations}
\begin{figure}[H]
  \centering
  \includegraphics[width=9cm]{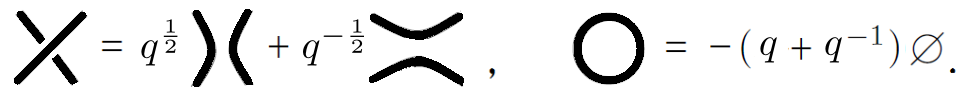}
\end{figure}
When $\Sigma$ is an oriented surface, $\mathcal{S}(\Sigma\times[0,1];R)$ is denoted by $\mathcal{S}(\Sigma;R)$. With $L_1L_2$ defined by stacking $L_1$ over $L_2$ in the $[0,1]$ direction, $\mathcal{S}(\Sigma;R)$ becomes a $R$-algebra, called the {\it Kauffman bracket skein algebra} of $\Sigma$.

When $R=\mathbb{C}$ and $q^{\frac{1}{2}}=-1$, by \cite{Bu97} there is a ring isomorphism
\begin{align}
\epsilon: \mathcal{S}(M;\mathbb{C})/{\rm nilradical}\to\mathbb{C}[\mathcal{X}(M)],  \label{eq:quantization}
\end{align}
where $\mathcal{X}(M)=\mathcal{X}(\pi_1(M))$ is the
${\rm SL}(2,\mathbb{C})$-character variety of $M$ (see Section \ref{sec:character-variety} for details).
Explicitly, $\epsilon$ sends a link $L=\sqcup_{i=1}^mK_i\subset M$ (the $K_i$'s being components) to the function
\begin{align}
\epsilon(L):\mathcal{X}(M)\to\mathbb{C}, \qquad \chi\mapsto{\prod}_{i=1}^m(-\chi([K_i])),  \label{eq:explicit}
\end{align}
where $[K_i]$ denotes the conjugacy class in $\pi_1(M)$ determined by $K_i$.

Let $\mathcal{S}(M)=\mathcal{S}(M;\mathbb{Z}[q^{\pm\frac{1}{2}}])$.
It was known \cite{BIMPW21,Pr99} that $\mathcal{S}(M)$ may admit torsion if $M$ contains an essential sphere or torus. As stated in
\cite[Problem 4.2]{Oh}, a question of fundamental importance is whether other surfaces can yield torsion as well.
Furthermore, the following problem remains open:
\begin{prob}[1.92 (G)-(i) in \cite{Ki97}]\label{prob:torsion}
Is $\mathcal{S}(M)$ free when $M$ is compact, irreducible and has no incompressible non-boundary parallel torus?
\end{prob}

\begin{figure}[H]
  \centering
  \includegraphics[width=3cm]{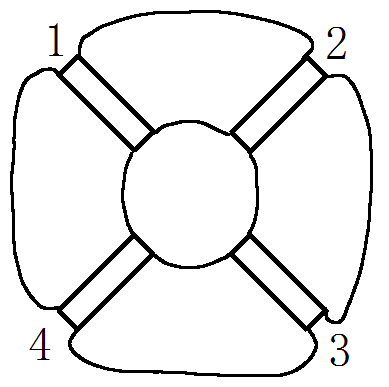} \\
  \caption{A $4$-strand Montesinos knot $K$. The rectangle labeled by $i$ stands for the rational tangle $[a_i/b_i]$,
  with $b_i\ge 3$.}\label{fig:Montesinos}
\end{figure}

Let $K=K(a_1/b_1,a_2/b_2,a_3/b_4,a_4/b_4)$, a $4$-strand Montesinos knot in $S^3$, as displayed in Figure \ref{fig:Montesinos}; we assume each $b_i\ge 3$.
Let $E_K$ denote the exterior of $K$, i.e. the complement of a tubular neighborhood of $K$ in $S^3$.

Let $e\in\mathcal{S}(E_K)$ denote the following element:
\begin{figure}[H]
  \centering
  \includegraphics[width=11.5cm]{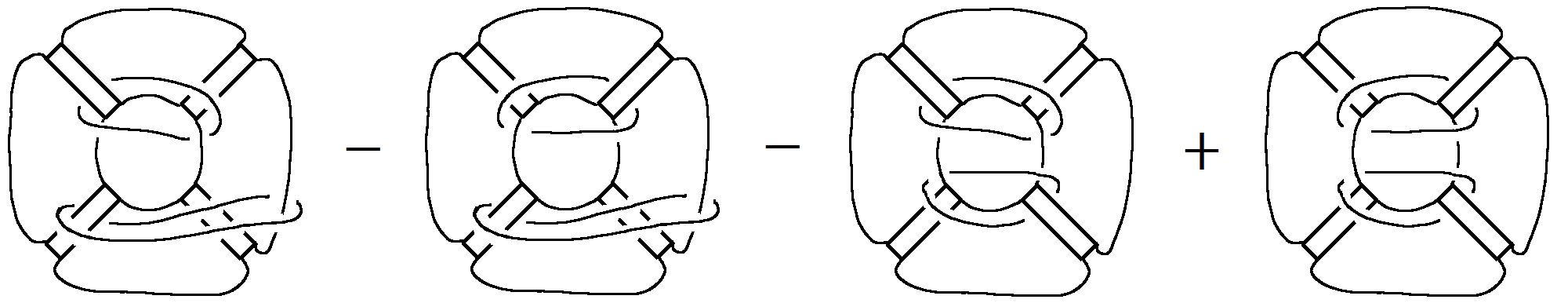}
\end{figure}

The main result of the paper is
\begin{thm}\label{thm:main}
The element $e$ does not vanish and satisfies $(q^2-1)e=0$.
\end{thm}

By \cite[Corollary 5]{Oe84}, $S^3\setminus K$ is hyperbolic, so $E_K$ provides a negative answer to Problem \ref{prob:torsion}.

Actually, incompressible surfaces in most Montesinos link exteriors were determined in \cite[Section 1]{Oe84}. By the results there, for the present knot $K$, except the peripheral torus and just one isotopy class which occurs when $\sum_{i=1}^4a_i/b_i=0$, each closed incompressible surface is isotopic to a surface obtained from a disjoint union of finitely many Conway spheres by a sequence of ``peripheral tubing" operations. Each Conway sphere is required to encircle exactly two tangles, and each tube is required to pass through at least one tangle. Then it follows that $E_K$ contains closed incompressible surface of every genus $\ge 2$, but does not contain an incompressible torus.
It is also known that each tubed incompressible surface in $E_K$ remains incompressible in any nontrivial Dehn filling of $E_K$.

Recently, Belletti and Detcherry \cite{BD25} gave criteria for the existence of torsion in $\mathcal{S}(M)$ for closed 3-manifolds $M$, based on the confirmation of Witten's finiteness conjecture, a highly nontrivial result of Gunningham, Jordan, and Safronov \cite{GJS23}. In contrast, our method of constructing torsion is straightforward and rather understandable.

The content is organized as follows.
In Section 2 we present many explicit computations in the skein module of the $3$-cube graph exterior, to obtain the fascinating identity (\ref{eq:wonderful}) which will play a key role, and show $(q^2-1)e=0$ in $\mathcal{S}(E_K)$.
To ensure the readability, we put several figures in the appendix. In Section 3, we show $e\ne 0$ so as to establish Theorem \ref{thm:main}, by studying a special class of ${\rm SL}(2,\mathbb{C})$-representations of $\pi_1(E_K)$. Moreover, we construct an infinite family of hyperbolic rational homology spheres whose skein modules admit torsion, by Dehn fillings of $E_K$.

Throughout the paper, we denote $q^{-1}$ as $\overline{q}$, denote $q^{-\frac{1}{2}}$ as $\overline{q}^{\frac{1}{2}}$, and so forth. Let $\alpha=q+\overline{q}$. Each link is equipped with the blackboard framing.

\medskip

\noindent
{\bf Acknowledgement}.

The author thanks Professor Thang T.Q. L\^e for valuable comments. He also thanks an anonymous referee for giving many constructive suggestions and comments, so that the paper can be largely improved.

\section{Skein calculations}\label{sec:graph}

\subsection{Set up}

\begin{figure}[h]
  \centering
  \includegraphics[width=9.5cm]{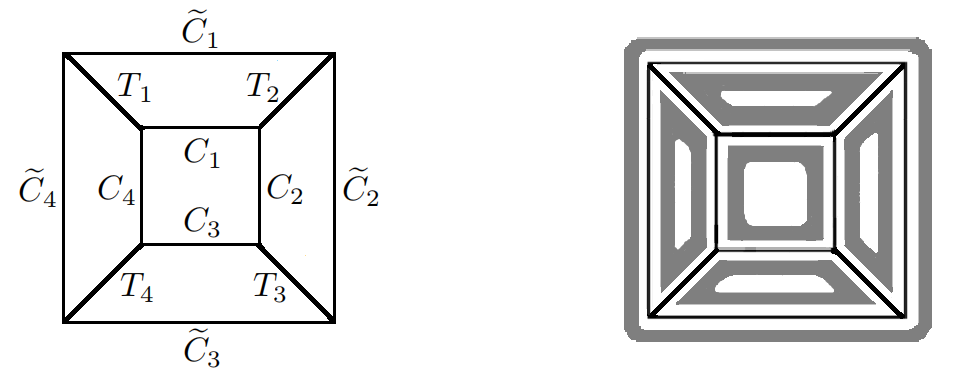} \\
  \caption{Left: the $3$-cube graph, whose edges are named as $C_i$, $\widetilde{C}_i$, $T_i$. We call the square made of $C_1,\ldots,C_4$ the ``inner square"; call the square made of $\widetilde{C}_1,\ldots,\widetilde{C}_4$ the ``outer square".
  Right: the shaded region displays a collar neighborhood $V$.}\label{fig:graph}
\end{figure}

Let $X$ denote the complement in $S^3$ of a tubular neighborhood of the $3$-cube graph shown in the left part of Figure \ref{fig:graph}.
Since $X$ is homeomorphic to a handlebody, by \cite[Theorem 2.3]{Pr99}, $\mathcal{S}(X)$ is free.

Let $\Sigma=\partial X$.
Take a closed collar neighborhood $V$ of $\Sigma$, as illustrated in the right part of Figure \ref{fig:graph}. Identify
$\Sigma\times[0,1]$ with $V$ in such a way that for $\mathsf{x}\in\Sigma,a\in[0,1]$, the distance from $(\mathsf{x},a)$ to $G$ grows as $a$ increases.
Recall that $\mathcal{S}(\Sigma)$ is a $R$-algebra.
A product on $\mathcal{S}(V)$ is induced through the identification $\mathcal{S}(V)\cong\mathcal{S}(\Sigma)$.

\begin{figure}[H]
  \centering
  \includegraphics[width=13cm]{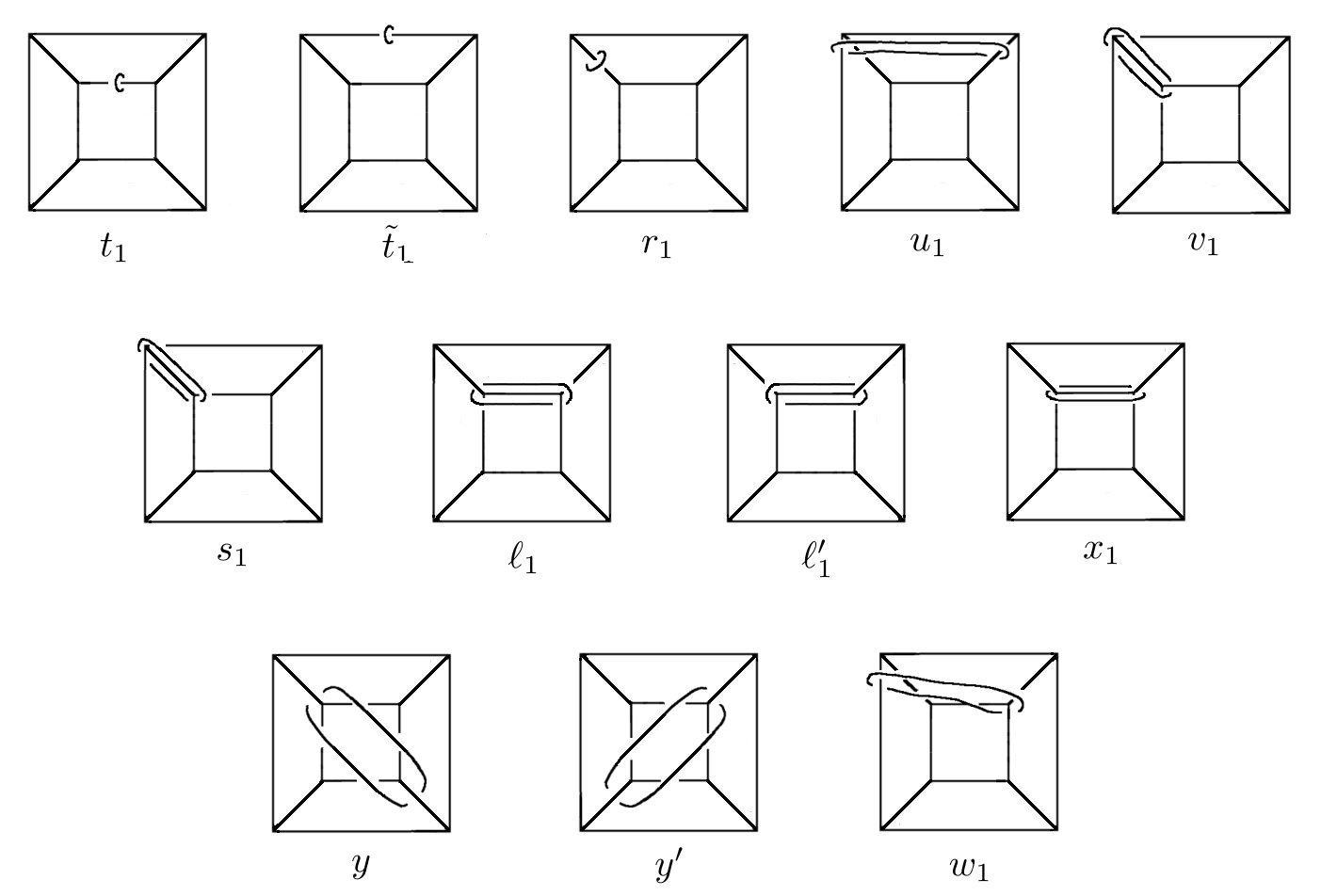} \\
  \caption{Upper and middle: some elements of $\mathcal{S}(V)$. Lower: some elements of $\mathcal{S}(X)$.}\label{fig:notation}
\end{figure}

Introduced in the upper and middle parts of Figure \ref{fig:notation} are some elements of $\mathcal{S}(V)$ which will play an important role.
For $2\le i\le 4$, let $t_i,\ldots,x_i$ respectively denote the elements obtained by rotating $t_1,\ldots,x_1$ by $(i-1)\pi/2$ counterclockwise.
We regard them as elements of $\mathcal{S}(X)$ via the morphism $\mathcal{S}(V)\to\mathcal{S}(X)$ induced by the inclusion $V\subset X$.
Note that the $r_i,t_i,\tilde{t}_i$ for $1\le i\le 4$ commute with each other, but for each $i$, none of $t_i,\ell_i,\ell'_i$ commutes with $x_i$.

Some other elements of $\mathcal{S}(X)$ are introduced in the lower part of Figure \ref{fig:notation}.
For $2\le i\le 4$, let $w_i$ denote the element obtained by rotating $w_1$ by $(i-1)\pi/2$ counterclockwise; it is a loop encircling $\widetilde{C}_{i-1}$ and $C_{i+1}$.

\begin{figure}[H]
  \centering
  \includegraphics[width=12.5cm]{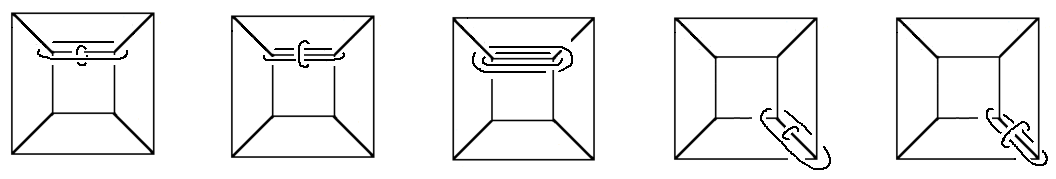} \\
  \caption{From left to right: $x_1t_1$, $t_1x_1$, $\ell_1x_1$, $s_3r_3$, $r_3s_3$.}\label{fig:example}
\end{figure}

Some examples of products in $\mathcal{S}(V)$ are given in Figure \ref{fig:example}.
Notice the difference between $x_1t_1$ and $t_1x_1$, as well as that between $s_3r_3$ and $r_3s_3$.

Take a homeomorphism ${\rm gl}:X\cup\Sigma\times[0,1]\cong X$. It makes $\mathcal{S}(X)$ into a right $\mathcal{S}(V)$-module.
Explicitly, for $b\in\mathcal{S}(X)$ and $a\in\mathcal{S}(V)\cong\mathcal{S}(\Sigma)$, one can define $ba\in\mathcal{S}(X)$ to be the image of $b\cup a\in\mathcal{S}(X\cup \Sigma\times[0,1])$ under the isomorphism induced by ${\rm gl}$.

\subsection{Basic formulas}

We shall frequently use the following local relations:
\begin{figure}[H]
  \centering
  \includegraphics[width=10.5cm]{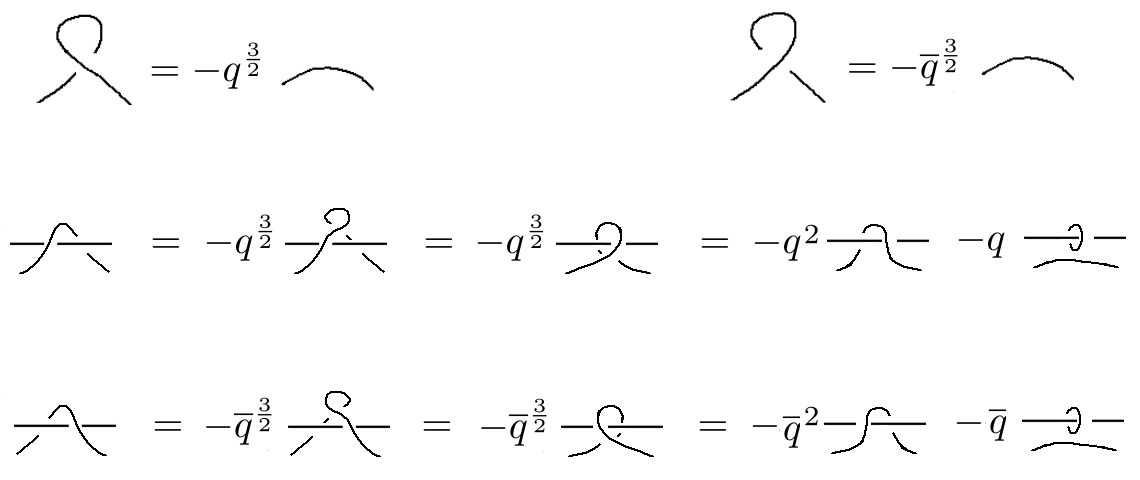}\\
  \caption{These are immediate consequences of the skein relations.}\label{fig:curl}
\end{figure}

\begin{figure}[H]
  \centering
  \includegraphics[width=13cm]{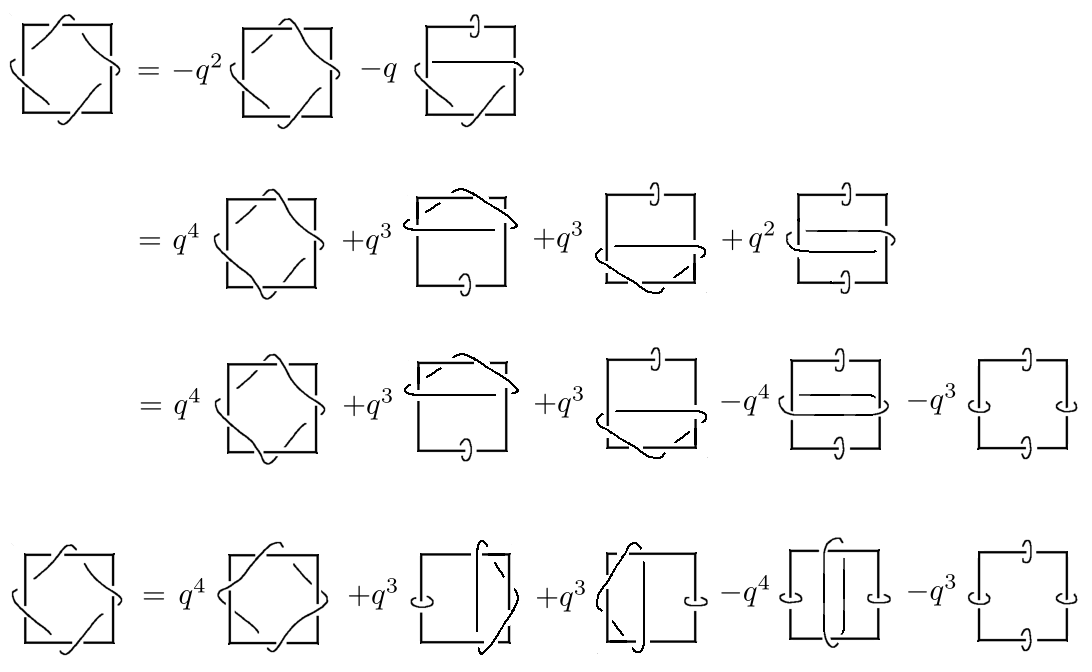}\\
  \caption{The upper formula is deduced by repeatedly applying the middle formula in Figure \ref{fig:curl}.
  The lower formula is obtained from rotating the upper one by $-\pi/2$.}\label{fig:difference-1}
\end{figure}

\begin{figure}[H]
  \centering
  \includegraphics[width=13cm]{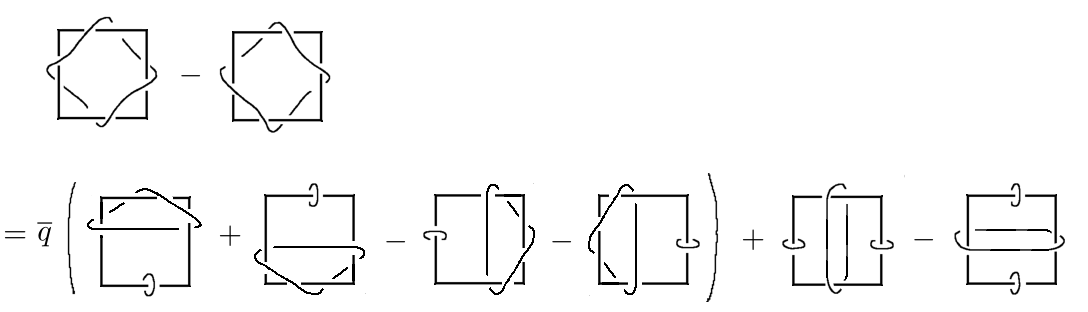}\\
  \caption{This is obtained by comparing the two formulas in Figure \ref{fig:difference-1}.}\label{fig:difference-2}
\end{figure}

\begin{figure}[H]
  \centering
  \includegraphics[width=13cm]{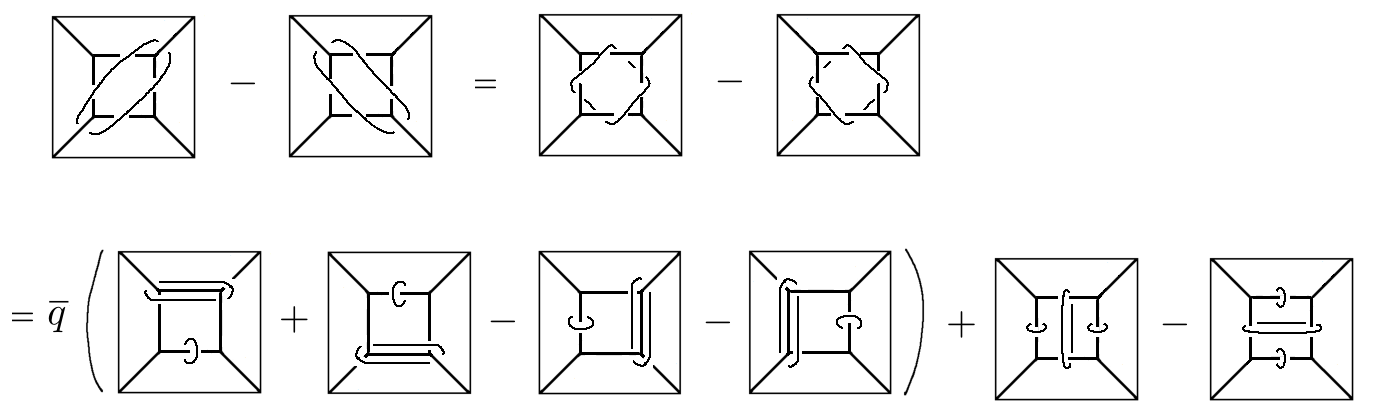}\\
  \caption{This results from applying the formula in Figure \ref{fig:difference-2} to the inner square.}\label{fig:y'-minus-y-1}
\end{figure}

\begin{figure}[H]
  \centering
  \includegraphics[width=13cm]{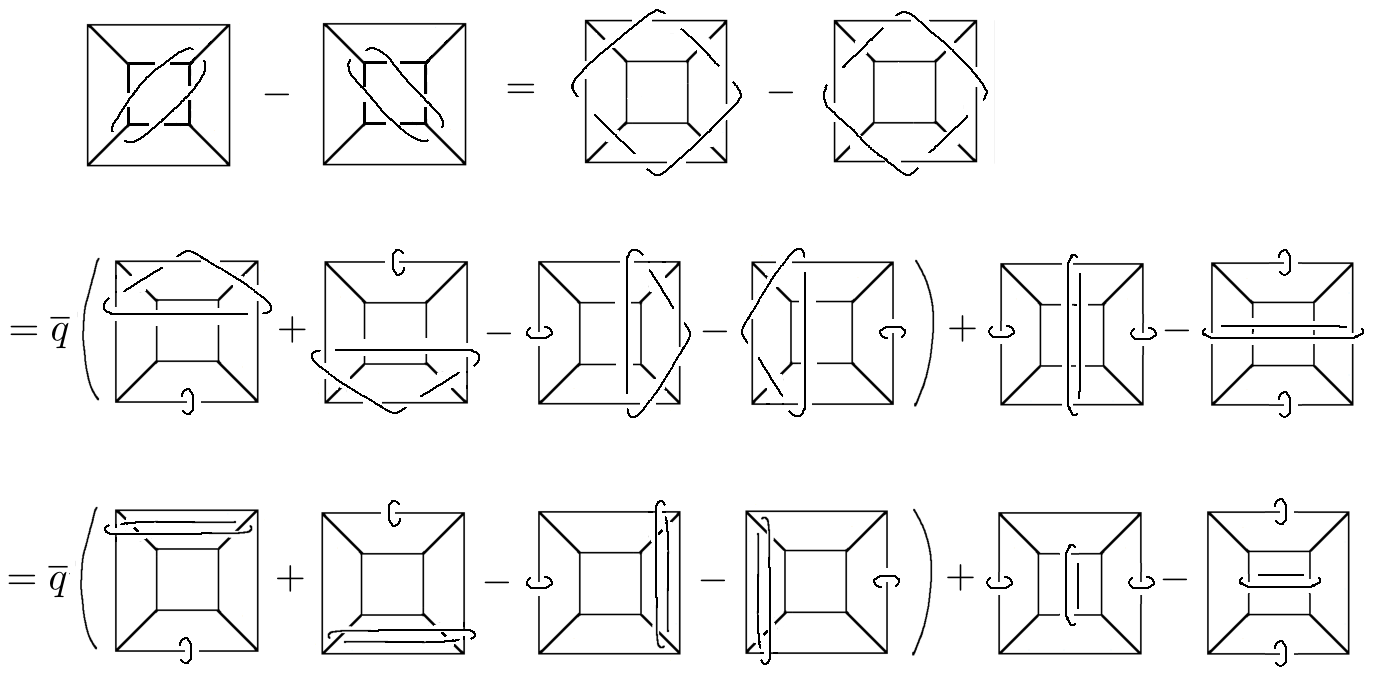}\\
  \caption{This results from applying the formula in Figure \ref{fig:difference-2} to the outer square.}\label{fig:y'-minus-y-2}
\end{figure}

From Figure \ref{fig:y'-minus-y-1}, Figure \ref{fig:y'-minus-y-2} we see respectively
\begin{align}
y'-y&=z:=\overline{q}(t_3\ell_1+t_1\ell_3-t_4\ell_2-t_2\ell_4)+x_2t_2t_4-x_1t_1t_3, \label{eq:y'-minus-y-1}  \\
y'-y&=\tilde{z}:=\overline{q}(\tilde{t}_4u_2+\tilde{t}_2u_4-\tilde{t}_3u_1-\tilde{t}_1u_3)+x_1\tilde{t}_1\tilde{t}_3-x_2\tilde{t}_2\tilde{t}_4.   \label{eq:y'-minus-y-2}
\end{align}

\subsection{A fascinating identity in $\mathcal{S}(X)$}

\begin{figure}[H]
  \centering
  \includegraphics[width=13cm]{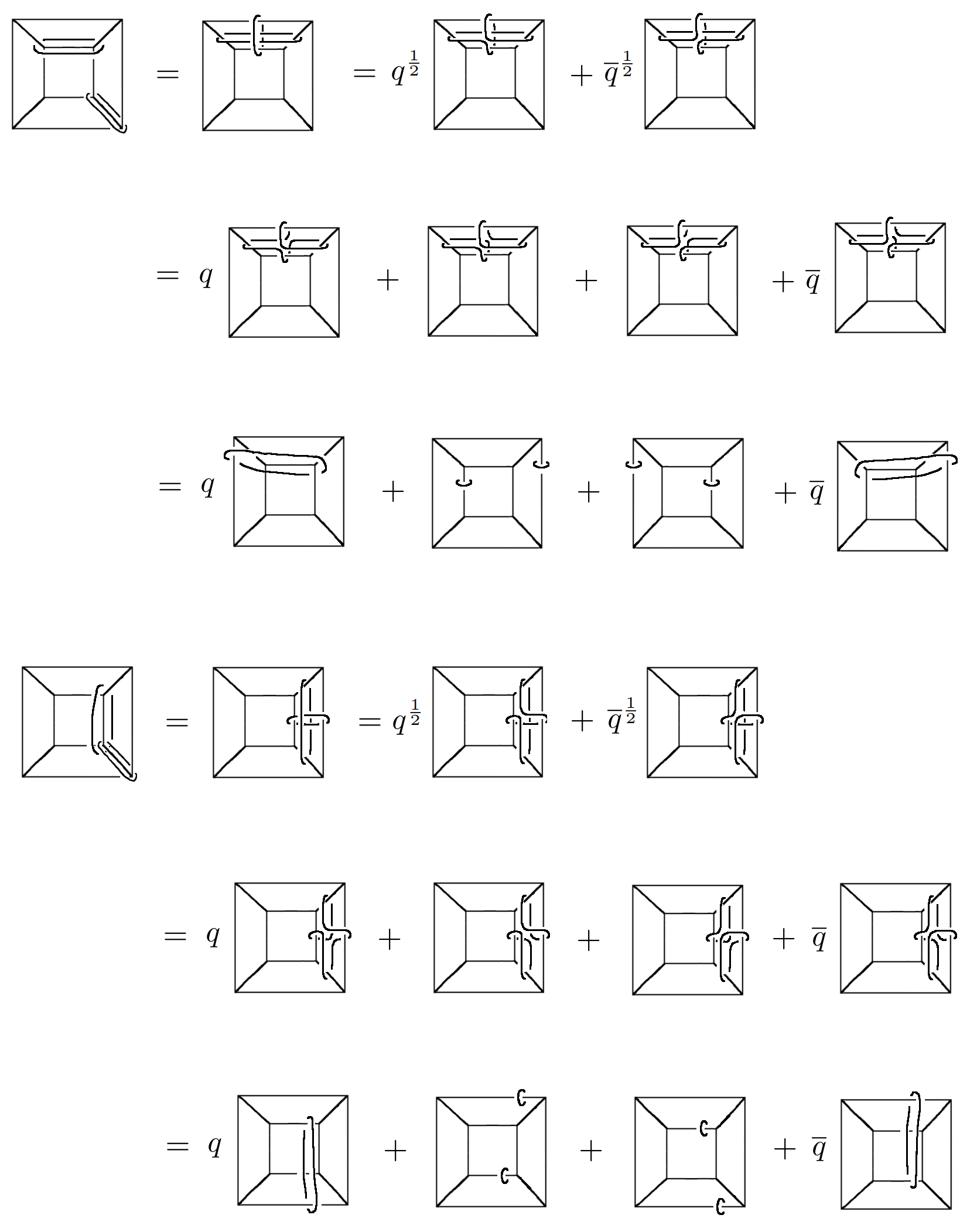}\\
  \caption{Computing $x_1s_3$ and $x_2s_3$.}\label{fig:isotopy-1}
\end{figure}

As illustrated in Figure \ref{fig:isotopy-1},
\begin{align*}
x_1s_3&=qw_1+\overline{q}w_3+\tilde{t}_2t_4+t_2\tilde{t}_4,   \\
x_2s_3&=qw_4+\overline{q}w_2+\tilde{t}_1t_3+t_1\tilde{t}_3.
\end{align*}

\begin{figure}[H]
  \centering
  \includegraphics[width=12.2cm]{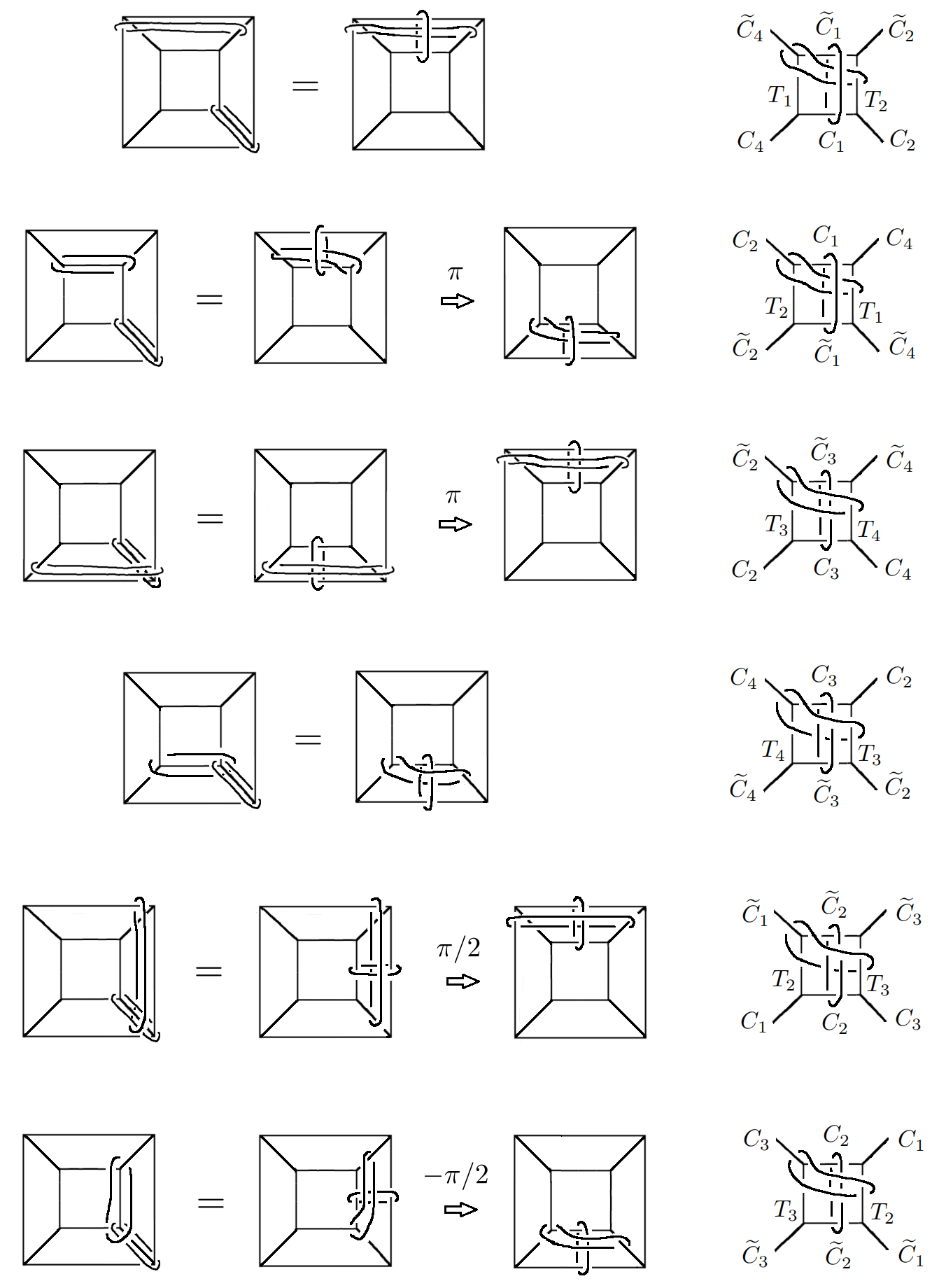}\\
  \caption{First to sixth rows: for computing $u_1s_3$, $\ell_1s_3$, $u_3s_3$, $\ell_3s_3$, $u_2s_3$, $\ell_2s_3$. The computations for $u_4s_3$, $\ell_4s_3$ are similar to $u_1s_3$, $\ell_1s_3$. In each row, displayed at the rightmost is the ``local shape" (just so named).}\label{fig:isotopy-2}
\end{figure}

\begin{figure}[H]
  \centering
  \includegraphics[width=12cm]{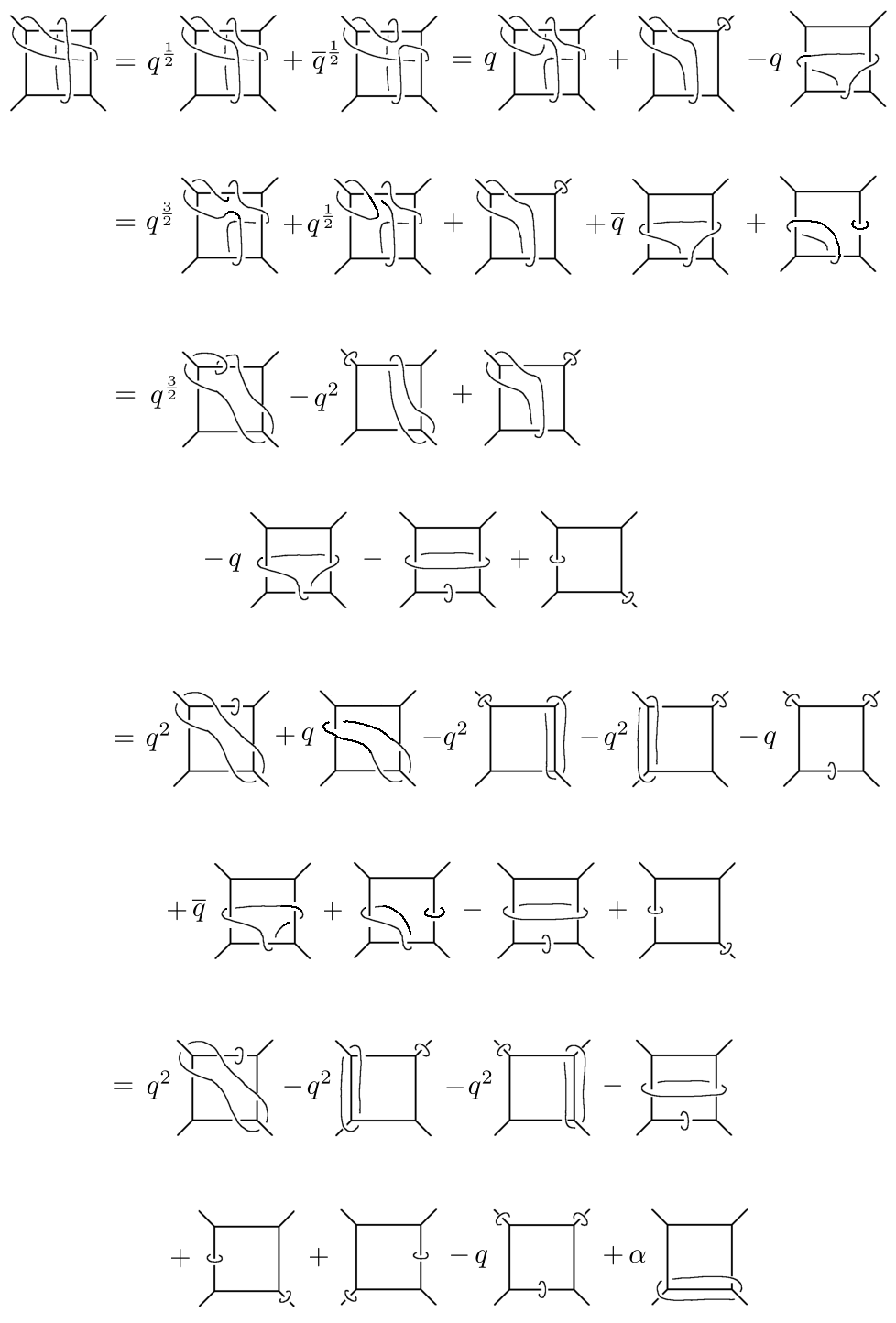}\\
  \caption{Applied to compute $u_1s_3$, $\ell_1s_3$, $u_4s_3$, $\ell_4s_3$.}\label{fig:expand-1}
\end{figure}

Applying the formula in Figure \ref{fig:expand-1} to the square made of $T_1,T_2,C_1,\widetilde{C}_1$, we obtain
\begin{align*}
u_1s_3&=q^2w_1\tilde{t}_1-q^2v_1\tilde{t}_2-q^2v_2\tilde{t}_4-x_1t_1+t_2r_1+t_4r_2-q\tilde{t}_4\tilde{t}_2t_1+\alpha\ell_1, \\
\ell_1s_3&=q^2w_1t_1-q^2v_2t_4-q^2v_1t_2-x_1\tilde{t}_1+\tilde{t}_4r_2+\tilde{t}_2r_1-qt_4t_2\tilde{t}_1+\alpha u_1.
\end{align*}
In the first equation, $s_3$ has been moved to encircle $u_1$; in the second equation, $s_3$ has been moved to encircle $\ell_1$, and we are actually applying the formula obtained from rotating the one in Figure \ref{fig:expand-1} by $\pi$.
Such conventions are also taken for other cases. See Figure \ref{fig:isotopy-2} for illustration.

\begin{figure}[H]
  \centering
  \includegraphics[width=12cm]{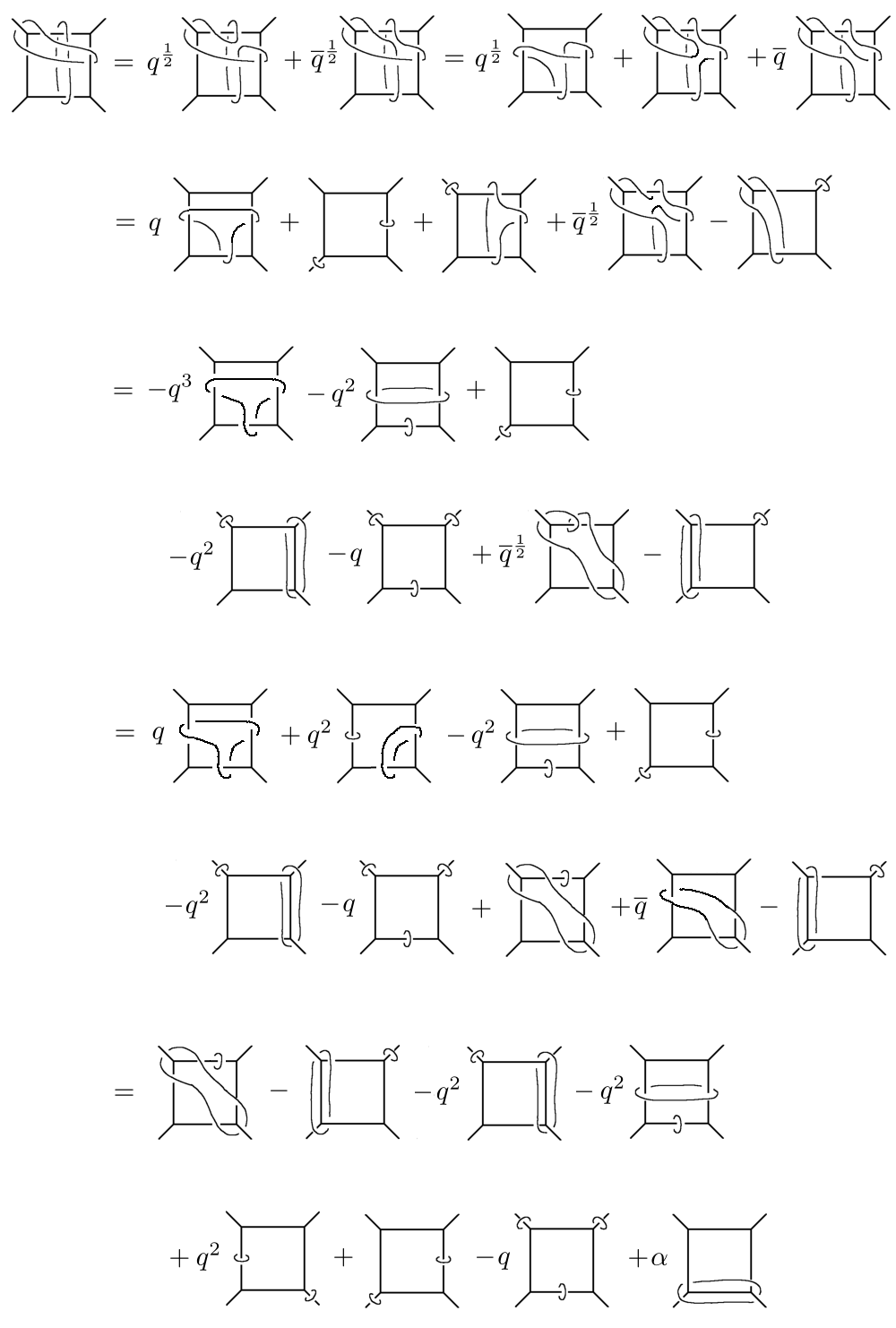}\\
  \caption{Applied to compute $u_3s_3$ and $\ell_2s_3$.}\label{fig:expand-2}
\end{figure}

\begin{figure}[H]
  \centering
  \includegraphics[width=12cm]{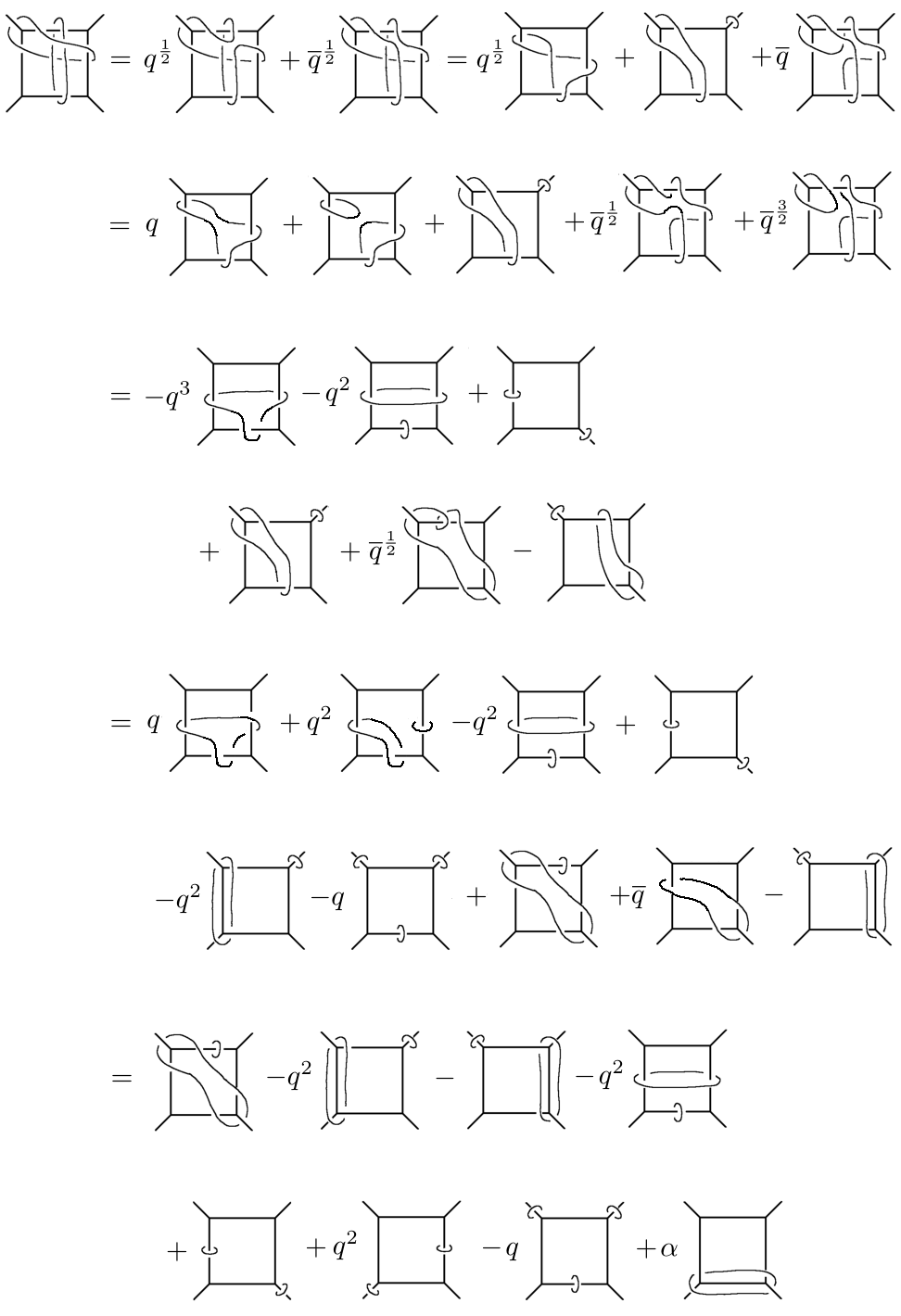}\\
  \caption{Applied to compute $\ell_3s_3$ and $u_2s_3$.}\label{fig:expand-3}
\end{figure}

By Figure \ref{fig:expand-2},
$$u_3s_3=w_3\tilde{t}_3-v_3\tilde{t}_4-q^2v_4\tilde{t}_2-q^2x_1t_3+q^2t_4r_3+t_2r_4-q\tilde{t}_2\tilde{t}_4t_3+\alpha\ell_3;$$
by Figure \ref{fig:expand-3},
$$\ell_3s_3=w_3t_3-q^2v_4t_2-v_3t_4-q^2x_1\tilde{t}_3+\tilde{t}_2r_4+q^2\tilde{t}_4r_3-qt_2t_4\tilde{t}_3+\alpha u_3.$$
Hence, remembering that $s_3$ commutes with $t_i,\tilde{t}_i$, we have
\begin{align*}
&\big(u_1\tilde{t}_3+\ell_1t_3+u_3\tilde{t}_1+\ell_3t_1-qx_1(t_1t_3+\tilde{t}_1\tilde{t}_3)\big)(s_3-\alpha) \\
=\ &\big(q^2w_1\tilde{t}_1-q^2v_1\tilde{t}_2-q^2v_2\tilde{t}_4-x_1t_1+t_2r_1+t_4r_2
-q\tilde{t}_4\tilde{t}_2t_1+\alpha\ell_1\big)\tilde{t}_3   \\
&+\big(q^2w_1t_1-q^2v_1t_2-q^2v_2t_4-x_1\tilde{t}_1+\tilde{t}_2r_1+\tilde{t}_4r_2-qt_4t_2\tilde{t}_1+\alpha u_1\big)t_3  \\
&+\big(w_3\tilde{t}_3-v_3\tilde{t}_4-q^2v_4\tilde{t}_2-q^2x_1t_3+q^2t_4r_3+t_2r_4
-q\tilde{t}_2\tilde{t}_4t_3+\alpha\ell_3\big)\tilde{t}_1  \\
&+\big(w_3t_3-v_3t_4-q^2v_4t_2-q^2x_1\tilde{t}_3+q^2\tilde{t}_4r_3+\tilde{t}_2r_4-qt_2t_4\tilde{t}_3+\alpha u_3\big)t_1  \\
&-q\big(qw_1+\overline{q}w_3+\tilde{t}_2t_4+t_2\tilde{t}_4\big)(t_1t_3+\tilde{t}_1\tilde{t}_3)   \\
&-\alpha(u_1\tilde{t}_3+\ell_1t_3+u_3\tilde{t}_1+\ell_3t_1)+q\alpha x_1(t_1t_3+\tilde{t}_1\tilde{t}_3) \\
=\ &-q^2v_1(\tilde{t}_2\tilde{t}_3+t_2t_3)-q^2v_2(\tilde{t}_3\tilde{t}_4+t_3t_4)
-v_3(\tilde{t}_1\tilde{t}_4+t_1t_4)-q^2v_4(\tilde{t}_1\tilde{t}_2+t_1t_2) \\
&+(\tilde{t}_2t_3+t_2\tilde{t}_3)r_1+(\tilde{t}_3t_4+t_3\tilde{t}_4)r_2+q^2(\tilde{t}_4t_1+t_4\tilde{t}_1)r_3+(\tilde{t}_1t_2+t_1\tilde{t}_2)r_4 \\
&-q(\tilde{t}_1t_3+t_1\tilde{t}_3)(\tilde{t}_2\tilde{t}_4+t_2t_4)-q(\tilde{t}_1\tilde{t}_3+t_1t_3)(\tilde{t}_2t_4+t_2\tilde{t}_4)  \\
&+q\alpha x_1(t_1-\tilde{t}_1)(t_3-\tilde{t}_3)+\alpha(u_1-\ell_1)(t_3-\tilde{t}_3)+\alpha(u_3-\ell_3)(t_1-\tilde{t}_1).
\end{align*}

By Figure \ref{fig:expand-3},
$$u_2s_3=w_2\tilde{t}_2-q^2v_2\tilde{t}_3-v_3\tilde{t}_1-q^2x_2t_2+t_3r_2+q^2t_1r_3-q\tilde{t}_1\tilde{t}_3t_2+\alpha\ell_2;$$ by Figure \ref{fig:expand-2},
$$\ell_2s_3=w_2t_2-v_3t_1-q^2v_2t_3-q^2x_2\tilde{t}_2+q^2\tilde{t}_1r_3+\tilde{t}_3r_2-qt_1t_3\tilde{t}_2+\alpha u_2;$$
by Figure \ref{fig:expand-1},
\begin{align*}
u_4s_3&=q^2w_4\tilde{t}_4-q^2v_4\tilde{t}_1-q^2v_1\tilde{t}_3-x_2t_4+t_1r_4+t_3r_1-q\tilde{t}_1\tilde{t}_3t_4+\alpha\ell_4, \\
\ell_4s_3&=q^2w_4t_4-q^2v_1t_3-q^2v_4t_1-x_2\tilde{t}_4+\tilde{t}_3r_1+\tilde{t}_1r_4-qt_1t_3\tilde{t}_4+\alpha u_4.
\end{align*}
Hence
\begin{align*}
&\big(u_2\tilde{t}_4+\ell_2t_4+u_4\tilde{t}_2+\ell_4t_2-qx_2(t_2t_4+\tilde{t}_2\tilde{t}_4)\big)(s_3-\alpha)  \\
=\ &\big(w_2\tilde{t}_2-q^2v_2\tilde{t}_3-v_3\tilde{t}_1-q^2x_2t_2+t_3r_2+q^2t_1r_3
-q\tilde{t}_1\tilde{t}_3t_2+\alpha\ell_2\big)\tilde{t}_4  \\
&+\big(w_2t_2-q^2v_2t_3-v_3t_1-q^2x_2\tilde{t}_2+\tilde{t}_3r_2+q^2\tilde{t}_1r_3-qt_1t_3\tilde{t}_2+\alpha u_2\big)t_4 \\
&+\big(q^2w_4\tilde{t}_4-q^2v_4\tilde{t}_1-q^2v_1\tilde{t}_3-x_2t_4+t_1r_4+t_3r_1
-q\tilde{t}_1\tilde{t}_3t_4+\alpha\ell_4\big)\tilde{t}_2  \\
&+\big(q^2w_4t_4-q^2v_4t_1-q^2v_1t_3-x_2\tilde{t}_4+\tilde{t}_1r_4+\tilde{t}_3r_1-qt_1t_3\tilde{t}_4+\alpha u_4\big)t_2  \\
&-q(qw_4+\overline{q}w_2+\tilde{t}_1t_3+t_1\tilde{t}_3)(t_2t_4+\tilde{t}_2\tilde{t}_4)  \\
&-\alpha(u_2\tilde{t}_4+\ell_2t_4+u_4\tilde{t}_2+\ell_4t_2)+q\alpha x_2(t_2t_4+\tilde{t}_2\tilde{t}_4)  \\
=\ &-q^2v_1(\tilde{t}_2\tilde{t}_3+t_2t_3)-q^2v_2(\tilde{t}_3\tilde{t}_4+t_3t_4)
-v_3(\tilde{t}_1\tilde{t}_4+t_1t_4)-q^2v_4(\tilde{t}_1\tilde{t}_2+t_1t_2)   \\
&+(\tilde{t}_2t_3+t_2\tilde{t}_3)r_1+(\tilde{t}_3t_4+t_3\tilde{t}_4)r_2+q^2(\tilde{t}_4t_1+t_4\tilde{t}_1)r_3+(\tilde{t}_1t_2+t_1\tilde{t}_2)r_4 \\
&-q(\tilde{t}_1\tilde{t}_3+t_1t_3)(\tilde{t}_2t_4+t_2\tilde{t}_4)-q(\tilde{t}_1t_3+t_1\tilde{t}_3)(\tilde{t}_2\tilde{t}_4+t_2t_4) \\
&+q\alpha x_2(t_2-\tilde{t}_2)(t_4-\tilde{t}_4)+\alpha(u_2-\ell_2)(t_4-\tilde{t}_4)+\alpha(u_4-\ell_4)(t_2-\tilde{t}_2).
\end{align*}

Therefore,
\begin{align*}
&q(z-\tilde{z})(s_3-\alpha)   \\
=\ &\big(u_1\tilde{t}_3+\ell_1t_3+u_3\tilde{t}_1+\ell_3t_1-qx_1(t_1t_3+\tilde{t}_1\tilde{t}_3))(s_3-\alpha)   \\
&-\big(u_2\tilde{t}_4+\ell_2t_4+u_4\tilde{t}_2+\ell_4t_2-qx_2(t_2t_4+\tilde{t}_2\tilde{t}_4)\big)(s_3-\alpha)   \\
=\ &\alpha\Big(qx_1(t_1-\tilde{t}_1)(t_3-\tilde{t}_3)-qx_2(t_2-\tilde{t}_2)(t_4-\tilde{t}_4)+(u_1-\ell_1)(t_3-\tilde{t}_3) \nonumber  \\
&\ \ \ \ +(u_3-\ell_3)(t_1-\tilde{t}_1)-(u_2-\ell_2)(t_4-\tilde{t}_4)-(u_4-\ell_4)(t_2-\tilde{t}_2)\Big).
\end{align*}
By (\ref{eq:y'-minus-y-1}), (\ref{eq:y'-minus-y-2}), $z=\tilde{z}$ in $\mathcal{S}(X)$, so the last line vanishes.
Since $\mathcal{S}(X)$ has no torsion, we have
\begin{align}
&qx_1(t_1-\tilde{t}_1)(t_3-\tilde{t}_3)-qx_2(t_2-\tilde{t}_2)(t_4-\tilde{t}_4)+(u_1-\ell_1)(t_3-\tilde{t}_3) \nonumber  \\
&+(u_3-\ell_3)(t_1-\tilde{t}_1)-(u_2-\ell_2)(t_4-\tilde{t}_4)-(u_4-\ell_4)(t_2-\tilde{t}_2)=0.  \label{eq:wonderful}
\end{align}

\subsection{Producing torsion in $\mathcal{S}(E_K)$}

\begin{figure}[H]
  \centering
  \includegraphics[width=12.2cm]{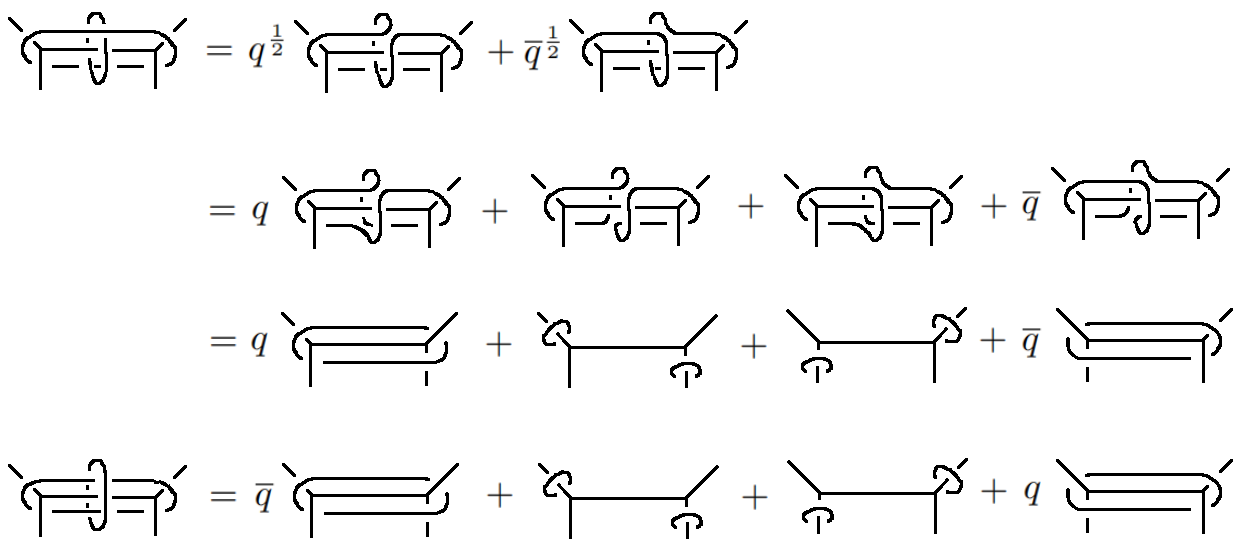}
  \caption{Upper: $x_1t_1=q\ell'_1+\overline{q}\ell_1+t_2r_1+t_4r_2$. Lower: $t_1x_1=\overline{q}\ell'_1+q\ell_1+t_2r_1+t_4r_2$, which can be deduced similarly as the upper one.}\label{fig:t-times-x}
\end{figure}

Remember that $x_1\tilde{t}_i=\tilde{t}_ix_1$ for all $i$ and $x_1t_i=t_ix_1$ for $i\ne 1$.
Multiplying (\ref{eq:wonderful}) by $x_1$ on the right yields
\begin{align}
&qx_1(t_1-\tilde{t}_1)x_1(t_3-\tilde{t}_3)-qx_2x_1(t_2-\tilde{t}_2)(t_4-\tilde{t}_4)+(u_1-\ell_1)x_1(t_3-\tilde{t}_3) \nonumber  \\
&+(u_3-\ell_3)x_1(t_1-\tilde{t}_1)-x_1(u_2-\ell_2)(t_4-\tilde{t}_4)-x_1(u_4-\ell_4)(t_2-\tilde{t}_2)    \nonumber  \\
&+(q-\overline{q})(u_3-\ell_3)(\ell_1-\ell'_1)=0;   \label{eq:near-torsion}
\end{align}
we have used the following identity in $\mathcal{S}(V)$:
\begin{align*}
t_1x_1-x_1t_1=(q-\overline{q})(\ell_1-\ell'_1),
\end{align*}
which can be obtained by comparing the two formulas in Figure \ref{fig:t-times-x}.

\begin{figure}[h]
  \centering
  \includegraphics[width=12cm]{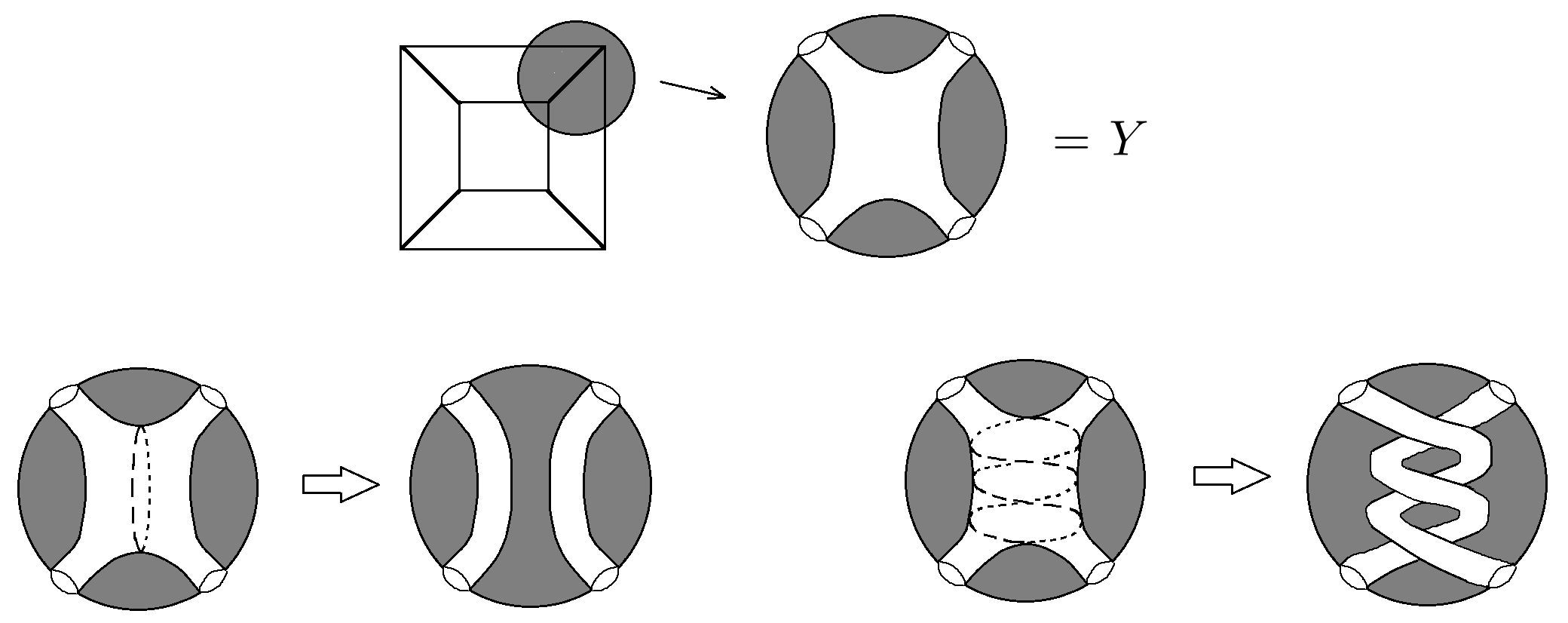}\\
  \caption{Upper: $Y=B\cap X$, where $B$ is a $3$-ball encircling a portion of the $3$-cube graph. Lower-left: attaching a $2$-handle to $Y$ along the dotted circle yields the exterior of a trivial tangle in $B$. Lower-right: considering the action via a certain self-homeomorphism of $Y$ preserving the boundary, we can see that attaching a $2$-handle to $Y$ along the dotted circle yields the exterior of a rational tangle in $B$.}\label{fig:explain}
\end{figure}

Now turn to $E_K$. Note that $E_K$ can be obtained by attaching four $2$-handles onto $X$, one for each rational tangle (see Figure \ref{fig:explain} for illustration).
By \cite[Proposition 2.2]{Pr99}, $\mathcal{S}(E_K)$ is a quotient of $\mathcal{S}(X)$,
so elements of $\mathcal{S}(X)$ can be used to represent elements of $\mathcal{S}(E_K)$.

A crucial ingredient is that when put at the innermost places (i.e. near $\partial E_K$), $\tilde{t}_i=t_i=t$ for all $i$, where $t$ is represented by a meridian of $K$.
This is caused by attaching $2$-handles, and this is why for each term in (\ref{eq:near-torsion}), we move at least one $t_i-\tilde{t}_i$ to the right end.

In $\mathcal{S}(E_K)$, (\ref{eq:near-torsion}) becomes
\begin{align}
(q^2-1)e=0, \qquad  e=(u_3-\ell_3)(\ell_1-\ell'_1).  \label{eq:torsion}
\end{align}

\begin{lem}\label{lem:torsion}
$e\ne 0$ in $\mathcal{S}(E_K)$.
\end{lem}

\section{Proving Lemma \ref{lem:torsion} using ${\rm SL}(2,\mathbb{C})$-representations}

\subsection{Preliminary}\label{sec:character-variety}

Let $\Gamma$ be a finitely presented group. A homomorphism $\rho:\Gamma\to{\rm SL}(2,\mathbb{C})$ is called a {\it representation}.
Call $\rho$ {\it reducible} if the elements in the image of $\rho$ share an eigenvector, and {\it irreducible} otherwise.
The {\it character} of $\rho$ is by definition the function $\chi_\rho:\Gamma\to\mathbb{C}$, $g\mapsto{\rm tr}(\rho(g))$.
The ${\rm SL}(2,\mathbb{C})$-{\it character variety} $\mathcal{X}(\Gamma)$ is defined as the GIT quotient
of $\hom(\Gamma,{\rm SL}(2,\mathbb{C}))$ by ${\rm SL}(2,\mathbb{C})$ acting via conjugation. As an affine algebraic set,
$$\mathcal{X}(\Gamma)=\{\chi_\rho\colon\rho\in\hom(\Gamma,{\rm SL}(2,\mathbb{C}))\}.$$
Let $\mathcal{X}^{\rm irr}(\Gamma)$ denote the set of characters of irreducible representations.
Call $\rho,\rho'$ {\it conjugate} if there exists $\mathbf{c}\in{\rm SL}(2,\mathbb{C})$ such that $\rho'(g)=\mathbf{c}\rho(g)\mathbf{c}^{-1}$ for all $g\in\Gamma$. Clearly, conjugate representations have the same character.

For a manifold $M$, denote $\mathcal{X}(M)=\mathcal{X}(\pi_1(M))$ and call it the character variety of $M$.
When $M$ is the exterior of a link $L$, also denote $\mathcal{X}(L)=\mathcal{X}(M)$.

Suppose $\mathsf{t}\in\mathbb{C}\setminus\{\pm2\}$. Let $G(\mathsf{t})=\{\mathbf{a}\in{\rm SL}(2,\mathbb{C})\colon{\rm tr}(\mathbf{a})=\mathsf{t}\}$.
Let
$$f_{\mathsf{t}}(\mathsf{r}_1,\mathsf{r}_2,\mathsf{r}_3;\mathsf{r})
=\mathsf{r}^2+\mathsf{t}(\mathsf{t}^2-\mathsf{r}_1-\mathsf{r}_2-\mathsf{r}_3)\mathsf{r}
+\mathsf{t}^2(3-\mathsf{r}_1-\mathsf{r}_2-\mathsf{r}_3)+\mathsf{r}_1^2+\mathsf{r}_2^2+\mathsf{r}_3^2+\mathsf{r}_1\mathsf{r}_2\mathsf{r}_3-4.$$

\begin{lem}\label{lem:matrix}
{\rm (i)} Two elements $\mathbf{a}_1,\mathbf{a}_2\in G(\mathsf{t})$ have a common eigenvector if and only if
${\rm tr}(\mathbf{a}_1\mathbf{a}_2)\in\{2,\mathsf{t}^2-2\}$.

{\rm (ii-1)} Given $\mathsf{t}_{12}\notin\{2,\mathsf{t}^2-2\}$, up to simultaneous conjugacy there exists a unique pair $(\mathbf{a}_1,\mathbf{a}_2)\in G(\mathsf{t})^{2}$ with ${\rm tr}(\mathbf{a}_1\mathbf{a}_2)=\mathsf{t}_{12}$.

{\rm (ii-2)} If furthermore $\mathsf{t}_{13},\mathsf{t}_{23},\mathsf{t}_{123}$ are given with $f_{\mathsf{t}}(\mathsf{t}_{12},\mathsf{t}_{13},\mathsf{t}_{23};\mathsf{t}_{123})=0$,
then there exists a unique $\mathbf{a}_3\in G(\mathsf{t})$ such that ${\rm tr}(\mathbf{a}_1\mathbf{a}_3)=\mathsf{t}_{13}$,
${\rm tr}(\mathbf{a}_2\mathbf{a}_3)=\mathsf{t}_{23}$ and ${\rm tr}(\mathbf{a}_1\mathbf{a}_2\mathbf{a}_3)=\mathsf{t}_{123}$.

{\rm(iii)} Given $\mathsf{t}_{41},\mathsf{t}_{12},\mathsf{t}_{23},\mathsf{t}_{34},\mathsf{t}_{124},\mathsf{t}_{234},\mathsf{t}_{24}\in\mathbb{C}$ with $\mathsf{t}_{24}\notin\{2,\mathsf{t}^2-2\}$ and
$$f_{\mathsf{t}}(\mathsf{t}_{41},\mathsf{t}_{12},\mathsf{t}_{24};\mathsf{t}_{124})
=f_{\mathsf{t}}(\mathsf{t}_{23},\mathsf{t}_{34},\mathsf{t}_{24};\mathsf{t}_{234})=0,$$
up to simultaneous conjugacy there exists a unique $4$-tuple $(\mathbf{x}_1,\ldots,\mathbf{x}_4)\in G(\mathsf{t})^{4}$ such that
${\rm tr}(\mathbf{x}_2\mathbf{x}_4)=\mathsf{t}_{24}$, ${\rm tr}(\mathbf{x}_1\mathbf{x}_2\mathbf{x}_4)=\mathsf{t}_{124}$,
${\rm tr}(\mathbf{x}_2\mathbf{x}_3\mathbf{x}_4)=\mathsf{t}_{234}$ and ${\rm tr}(\mathbf{x}_{i-1}\mathbf{x}_i)=\mathsf{t}_{i-1,i}$ for all $i$.
\end{lem}

\begin{proof}
For (i), (ii-1), (ii-2), see \cite{Go09}, in particular, Section 2 and Section 5.

(iii) By (ii-1), up to simultaneous conjugacy there exists a unique pair
$(\mathbf{x}_2,\mathbf{x}_4)\in G(\mathsf{t})^{2}$ with ${\rm tr}(\mathbf{x}_2\mathbf{x}_4)=\mathsf{t}_{24}$.
Then by (ii-2), there exists a unique $\mathbf{x}_1\in G(\mathsf{t})$ such that ${\rm tr}(\mathbf{x}_4\mathbf{x}_1)=\mathsf{t}_{41}$,
${\rm tr}(\mathbf{x}_1\mathbf{x}_2)=\mathsf{t}_{12}$ and ${\rm tr}(\mathbf{x}_1\mathbf{x}_2\mathbf{x}_4)=\mathsf{t}_{124}$; also, there exists a unique $\mathbf{x}_3\in G(\mathsf{t})$ such that ${\rm tr}(\mathbf{x}_2\mathbf{x}_3)=\mathsf{t}_{23}$,
${\rm tr}(\mathbf{x}_3\mathbf{x}_4)=\mathsf{t}_{34}$ and ${\rm tr}(\mathbf{x}_2\mathbf{x}_3\mathbf{x}_4)=\mathsf{t}_{234}$.
\end{proof}

Refer to \cite[Section 2]{FP16} for a nice exposition on rational tangles and Montesinos links.
To each $a/b\in\mathbb{Q}$ with $(a,b)=1$ is associated a rational tangle $[a/b]$ obtained by drawing slope $a/b$ lines on a ``square pillowcase" starting at the four corners (labeled by NW, NE, SW, SE). Using the expression of $a/b$ as a continued fraction, there is an alternative presentation for $[a/b]$ which is sometimes more convenient.
Connecting NW to SW and NE to SE by arcs, one obtains a rational link $B(a/b)$. When $b$ is odd, $B(a/b)$ is a knot; otherwise, $B(a/b)$ has two components.

\begin{figure}[H]
  \centering
  \includegraphics[width=9cm]{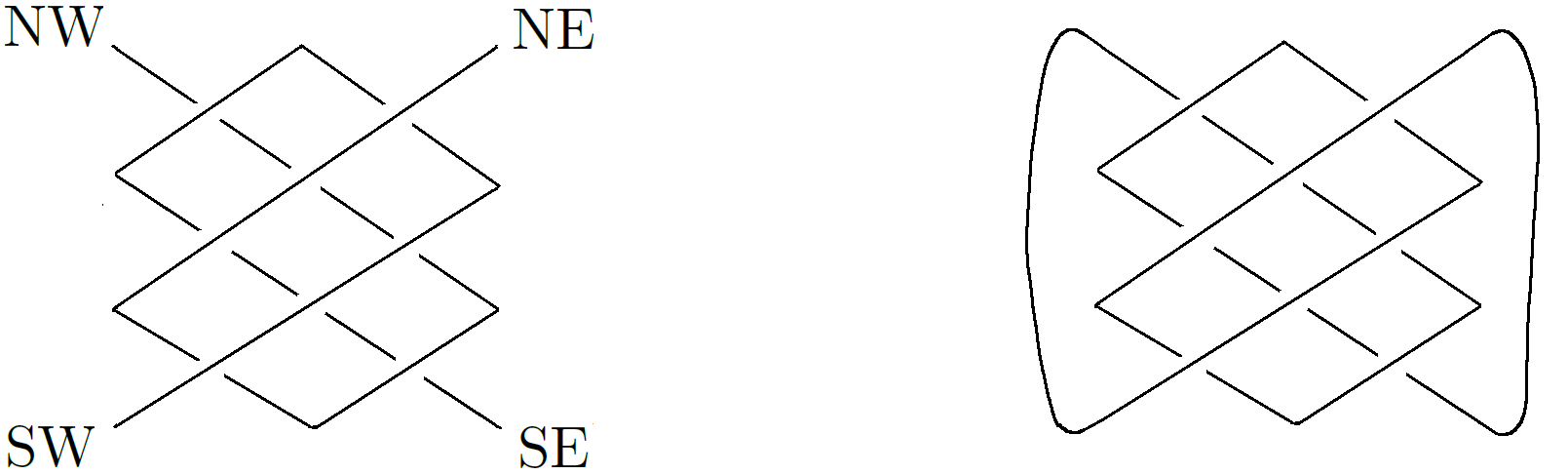}\\
  \caption{Left: the rational tangle [2/3]. Right: the rational knot $B(2/3)$.}\label{fig:tangle}
\end{figure}

A Montesinos link is the cyclic sum of a finite ordered list of rational tangles.
In particular, our knot $K=K(a_1/b_1,a_2/b_2,a_3/b_3,a_4/b_4)$
is constructed by joining the NE (resp. SE) of $[a_i/b_i]$ to the NW (resp. SW) of $[a_{i+1}/b_{i+1}]$ for $i=1,2,3,4$.

\begin{lem}\label{lem:rep}
Let $b\ge 3$. 
Given a generic $\mathsf{t}$, the rational link $B(a/b)$ admits an irreducible representation sending each meridian to an element of $G(\mathsf{t})$.
\end{lem}
One may refer to \cite[Section 2]{Ri84}, \cite[Section 2]{Ri92} and \cite[Section 7]{ORS08}.

Understood alternatively, it is known that $B(a/b)$ is hyperbolic, so by Thurston's result, the irreducible component of $\mathcal{X}(B(a/b))$ containing the character of a lift of the holonomy representation has dimension $1$ (resp. $2$) when $b$ is odd (resp. even).

\begin{rmk}
\rm Throughout this section, ``generic" means that finitely many values are excluded. In Lemma \ref{lem:rep}, those are given by $\mathsf{t}=\kappa+\kappa^{-1}$, where $\kappa^2$ are the roots of the Alexander polynomial of $B(a/b)$, corresponding to nonabelian but reducible representations of $B(a/b)$.
\end{rmk}

For a link $L$, call a representation $\pi_1(S^3\setminus L)\to{\rm SL}(2,\mathbb{C})$ simply a representation of $L$.
Given a projection diagram of $L$, let $\mathcal{D}(L)$ denote the set of directed arcs; each arc gives rise to two directed arcs $\mathfrak{a},\mathfrak{a}^{-1}$. From the viewpoint of Wirtinger presentation, a representation of $L$ can be identified with a map $\rho:\mathcal{D}(L)\to{\rm SL}(2,\mathbb{C})$ such that $\rho(\mathfrak{a}^{-1})=\rho(\mathfrak{a})^{-1}$ for each $\mathfrak{a}$ and $\rho(\mathfrak{c})=\rho(\mathfrak{a})\rho(\mathfrak{b})\rho(\mathfrak{a})^{-1}$ for all $\mathfrak{a},\mathfrak{b},\mathfrak{c}$
forming a crossing as follows:
\begin{figure}[H]
  \centering
  \includegraphics[width=2.8cm]{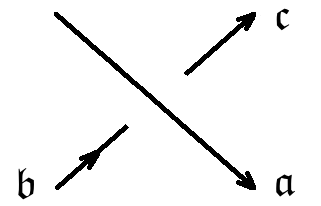}
\end{figure}

By this manner, we can also deal with representations of tangles.

\subsection{A subvariety of $\mathcal{X}(E_K)$}

Let $\mathcal{X}^\ast$ denote the subset of $\mathcal{X}(K)$ consisting of the characters of representations $\rho$ of $K$ such that $\rho(\tilde{\mathfrak{x}}_1)\rho(\mathfrak{x}_1)=\mathbf{e}$, where $\mathbf{e}$ is the identity matrix.

Suppose a generic $\mathsf{t}\ne\pm 2$ is given.
For each $i$, by Lemma \ref{lem:rep} we may take an irreducible representation $\varrho_i$ of $B(a_i/b_i)$ sending each meridian into $G(\mathsf{t})$; it can be identified with a representation of $[a_i/b_i]$, which is by abuse of notation also denoted by $\varrho_i$, such that $\varrho_i(\tilde{\mathfrak{x}}_{i-1})=\varrho_i(\mathfrak{x}_{i-1})^{-1}$ and $\varrho_i(\tilde{\mathfrak{x}}_i)=\varrho_i(\mathfrak{x}_i)^{-1}$. Let $\mathbf{u}_i=\varrho_i(\mathfrak{x}_{i-1})$, $\mathbf{v}_i=\varrho_i(\mathfrak{x}_{i})$, and let
$\mathsf{s}_i={\rm tr}(\mathbf{u}_i\mathbf{v}_{i})$. Note that $\mathbf{u}_i,\mathbf{v}_i$ have no common eigenvector; otherwise $\varrho_i$ would be reducible, since by \cite[Lemma 3.1]{Ch22}, each element in the image of $\varrho_i$ is a word in $\mathbf{u}_i$, $\mathbf{v}_i$. Then by Lemma \ref{lem:matrix} (i), $\mathsf{s}_i\notin\{2,\mathsf{t}^2-2\}$.

\begin{figure}[h]
  \centering
  \includegraphics[width=4.5cm]{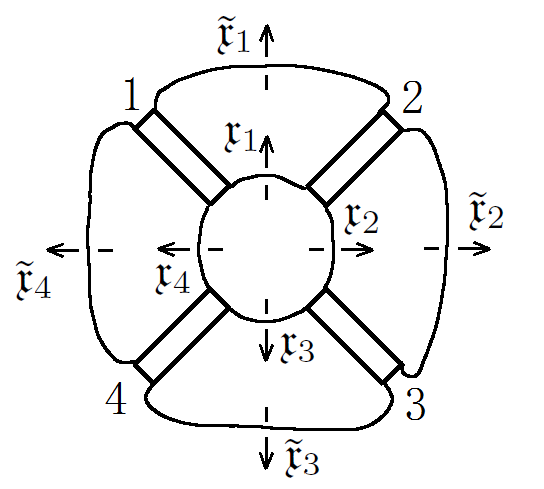}\\
  \caption{In the Wirtinger presentation, the $\mathfrak{x}_i$'s and $\tilde{\mathfrak{x}}_i$'s denote certain elements of $\pi_1(E_K)$.}\label{fig:representation}
\end{figure}

From these $\varrho_i$'s together with $\mathsf{c}\notin\{2,\mathsf{t}^2-2\}$, we can cook up a representation $\rho^{\mathsf{t},\mathsf{c}}$ of $K$.
By Lemma \ref{lem:matrix} (iii), there exist $\mathbf{x}_1,\ldots,\mathbf{x}_4\in G(\mathsf{t})$ with
${\rm tr}(\mathbf{x}_2\mathbf{x}_4)=\mathsf{c}$ and ${\rm tr}(\mathbf{x}_{i-1}\mathbf{x}_i)=\mathsf{s}_i$ for all $i$.
For each $i$, since ${\rm tr}(\mathbf{u}_i\mathbf{v}_i)={\rm tr}(\mathbf{x}_{i-1}\mathbf{x}_i)$, by Lemma \ref{lem:matrix} (ii-1) there exists $\mathbf{c}_i\in{\rm SL}(2,\mathbb{C})$ such that
$\mathbf{c}_i\mathbf{u}_i\mathbf{c}_i^{-1}=\mathbf{x}_{i-1}$ and $\mathbf{c}_i\mathbf{v}_i\mathbf{c}_i^{-1}=\mathbf{x}_i$.
Let $\rho_i=\mathbf{c}_i\varrho_i\mathbf{c}_i^{-1}$, by which we mean the representation of $[a_i/b_i]$ sending each directed arc $\mathfrak{a}$ to $\mathbf{c}_i\varrho_i(\mathfrak{a})\mathbf{c}_i^{-1}$. The $\rho_i$'s fit together to define the desired representation $\rho^{\mathsf{t},\mathsf{c}}$, whose character belongs to $\mathcal{X}^\ast$.

Thus, $\mathcal{X}^\ast$ is an algebraic surface, locally parameterized by $\mathsf{t},\mathsf{c}$.

Elements of $\mathcal{X}^\ast$ are identified with simultaneous conjugacy classes of $(\mathbf{x}_1,\ldots,\mathbf{x}_4)$'s subject to certain conditions.
For $i_1,\ldots,i_s\in\{1,\ldots,4\}$, let $\mathsf{t}_{i_1\cdots i_s}\in\mathbb{C}[\mathcal{X}^\ast]$ denote
the function sending $(\mathbf{x}_1,\ldots,\mathbf{x}_4)$ to ${\rm tr}(\mathbf{x}_{i_1}\cdots\mathbf{x}_{i_s})$.

\subsection{Proof of Lemma \ref{lem:torsion} and discussions}\label{sec:proof}

Fix a generic $\mathsf{t}\ne\pm2$, and fix $\mathsf{s}_1,\ldots,\mathsf{s}_4$ as above. Let $\mathcal{X}_1\subset\mathcal{X}^\ast$ denote the Zariski closure of the subset $\{\chi_{\rho^{\mathsf{t},\mathsf{c}}}\colon \mathsf{c}\ne 2,\mathsf{t}^2-2\}$.
Note that on $\mathcal{X}_1$, we have $\mathsf{t}_i\equiv\mathsf{t}$ and $\mathsf{t}_{i-1,i}\equiv\mathsf{t}^2-\mathsf{s}_i$ for each $i$, while $\mathsf{t}_{24}$ is nonconstant.

\begin{proof}[Proof of Lemma \ref{lem:torsion}]
Set $q^{\frac{1}{2}}=-1$. Recall (\ref{eq:quantization}), (\ref{eq:explicit}) that through $\epsilon$, elements of $\mathcal{S}(E_K)$ can be regarded as functions on $\mathcal{X}(E_K)$, so they can be regarded as functions on $\mathcal{X}_1$.
In particular, $\epsilon(t)=-\mathsf{t}$.

\begin{figure}[h]
  \centering
  \includegraphics[width=12.5cm]{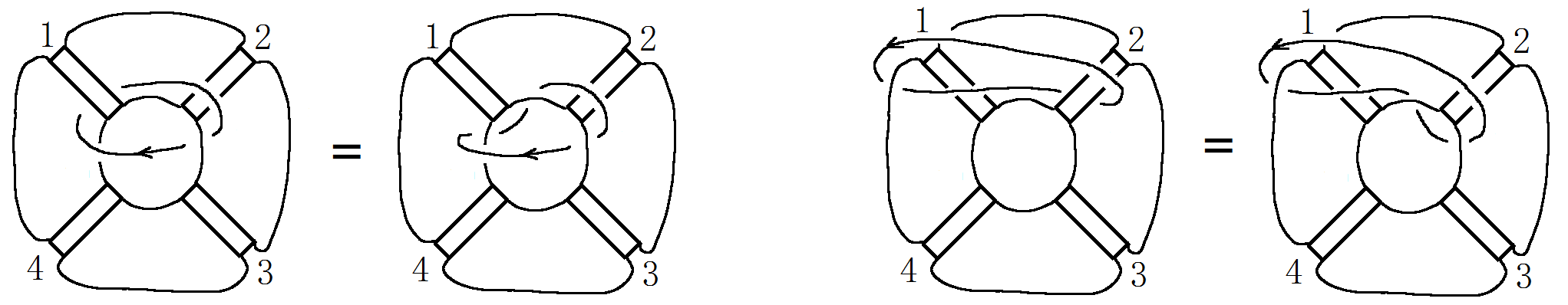}\\
  \caption{Arbitrarily choose orientations for $\ell_1,u_1$. Given $(\mathbf{x}_1,\ldots,\mathbf{x}_4)$ which determines a representation $\rho$ with $\chi_\rho\in\mathcal{X}^\ast$, the values of $\epsilon(\ell_1)$ and $\epsilon(u_1)$ taken at $(\mathbf{x}_1,\ldots,\mathbf{x}_4)$ are
  $-{\rm tr}(\rho(\mathfrak{x}_4^{-1}\mathfrak{x}_1\mathfrak{x}_2^{-1}))=-{\rm tr}(\mathbf{x}_4^{-1}\mathbf{x}_1\mathbf{x}_2^{-1})$ and
  $-{\rm tr}(\rho(\tilde{\mathfrak{x}}_4^{-1}\mathfrak{x}_1^{-1}\mathfrak{x}_2))
  =-{\rm tr}(\mathbf{x}_4\mathbf{x}_1^{-1}\mathbf{x}_2)$, respectively. Remember that $\rho(\tilde{\mathfrak{x}}_4)=\mathbf{x}_4^{-1}$.}\label{fig:map}
\end{figure}

Using $\mathbf{x}_i^{-1}=\mathsf{t}\mathbf{e}-\mathbf{x}_i$, we can compute that
$\epsilon(\ell_i)$ sends $(\mathbf{x}_1,\ldots,\mathbf{x}_4)$ to
\begin{align*}
-{\rm tr}\big(\mathbf{x}_{i-1}^{-1}\mathbf{x}_i\mathbf{x}_{i+1}^{-1}\big)
&={\rm tr}\big((\mathbf{x}_{i-1}-\mathsf{t}\mathbf{e})\mathbf{x}_i(\mathsf{t}\mathbf{e}-\mathbf{x}_{i+1})\big)   \\
&=\mathsf{t}\cdot{\rm tr}(\mathbf{x}_{i-1}\mathbf{x}_i)-{\rm tr}(\mathbf{x}_{i-1}\mathbf{x}_i\mathbf{x}_{i+1})
-\mathsf{t}^2\cdot{\rm tr}(\mathbf{x}_i)+\mathsf{t}\cdot{\rm tr}(\mathbf{x}_i\mathbf{x}_{i+1}),
\end{align*}
and $\epsilon(\ell'_i)$, $\epsilon(u_i)$ both send $(\mathbf{x}_1,\ldots,\mathbf{x}_4)$ to
$$-{\rm tr}\big(\mathbf{x}_{i-1}\mathbf{x}_i^{-1}\mathbf{x}_{i+1}\big)
={\rm tr}\big(\mathbf{x}_{i-1}(\mathbf{x}_i-\mathsf{t}\mathbf{e})\mathbf{x}_{i+1}\big)
={\rm tr}(\mathbf{x}_{i-1}\mathbf{x}_i\mathbf{x}_{i+1})-\mathsf{t}\cdot{\rm tr}(\mathbf{x}_{i-1}\mathbf{x}_{i+1}).$$
See Figure \ref{fig:map} for the illustration of $\epsilon(\ell_1)$, $\epsilon(u_1)$. The other cases are similar.

Hence
\begin{align*}
\epsilon(\ell_i)&=\mathsf{t}(\mathsf{t}_{i,i+1}+\mathsf{t}_{i-1,i}-\mathsf{t}^2)-\mathsf{t}_{i-1,i,i+1},    \\
\epsilon(\ell'_i)=\epsilon(u_i)&=\mathsf{t}_{i-1,i,i+1}-\mathsf{t}\mathsf{t}_{i-1,i+1}.
\end{align*}
Consequently,
\begin{align*}
\epsilon(\ell'_1-\ell_1)&=2\mathsf{t}_{124}-\mathsf{t}(\mathsf{t}_{12}+\mathsf{t}_{41}+\mathsf{t}_{24}-\mathsf{t}^2),   \\
\epsilon(u_3-\ell_3)&=2\mathsf{t}_{234}-\mathsf{t}(\mathsf{t}_{23}+\mathsf{t}_{34}+\mathsf{t}_{24}-\mathsf{t}^2).
\end{align*}

If $\epsilon(\ell'_1-\ell_1)=0$, then $\mathsf{t}_{124}=(1/2)\mathsf{t}(\mathsf{t}_{41}+\mathsf{t}_{12}+\mathsf{t}_{24}-\mathsf{t}^2)$, which combined with
$f_{\mathsf{t}}(\mathsf{t}_{41},\mathsf{t}_{12},\mathsf{t}_{24};\mathsf{t}_{124})=0$ implies 
\begin{align}
&\Big(\frac{\mathsf{t}^2}{4}-1\Big)\mathsf{t}_{24}^2+\Big(\frac{\mathsf{t}^2}{2}(\mathsf{t}_{41}+\mathsf{t}_{12}-\mathsf{t}^2)
+\mathsf{t}^2-\mathsf{t}_{41}\mathsf{t}_{12}\Big)\mathsf{t}_{24} \nonumber \\
&\ \ \ +\frac{\mathsf{t}^2}{4}(\mathsf{t}_{41}+\mathsf{t}_{12}-\mathsf{t}^2)^2+\mathsf{t}^2(\mathsf{t}_{41}+\mathsf{t}_{12}-3)
-(\mathsf{t}_{41}^2+\mathsf{t}_{12}^2-4)=0.  \label{eq:quadratic-1}
\end{align}
Hence $\epsilon(\ell'_1-\ell_1)\ne0$ in a Zariski open subset of $\mathcal{X}_1$.

Similarly, $\epsilon(u_3-\ell_3)=0$ will lead to
\begin{align}
&\Big(\frac{\mathsf{t}^2}{4}-1\Big)\mathsf{t}_{24}^2+\Big(\frac{\mathsf{t}^2}{2}(\mathsf{t}_{23}+\mathsf{t}_{34}-\mathsf{t}^2)
+\mathsf{t}^2-\mathsf{t}_{23}\mathsf{t}_{34}\Big)\mathsf{t}_{24} \nonumber \\
&\ \ \ +\frac{\mathsf{t}^2}{4}(\mathsf{t}_{23}+\mathsf{t}_{34}-\mathsf{t}^2)^2+\mathsf{t}^2(\mathsf{t}_{23}+\mathsf{t}_{34}-3)
-(\mathsf{t}_{23}^2+\mathsf{t}_{34}^2-4)=0.  \label{eq:quadratic-2}
\end{align}
Hence $\epsilon(u_3-\ell_3)\ne0$ in a Zariski open subset of $\mathcal{X}_1$.

Therefore, $\epsilon(e)\ne0$ in $\mathbb{C}[\mathcal{X}_1]$, implying $e\ne0$ in $\mathcal{S}(E_K)$.
\end{proof}

\begin{rmk}
\rm When $b_i\le 2$ for some $i$, the identities (\ref{eq:wonderful}) and (\ref{eq:torsion}) still hold, but we are no longer able to find a subvariety of $\mathcal{X}(K)$ such as $\mathcal{X}_1$ to rule out $\epsilon(e)=0$, so $e=0$ is probable.
\end{rmk}

As a negative answer to Problem 1.1 in the case of closed $3$-manifolds, \cite[Corollary 3.1]{BD25} asserts that the skein module of any closed hyperbolic $3$-manifold with positive first Betti number admits torsion. In comparison, we are able to construct hyperbolic rational homology spheres whose skein modules admit torsion.

Take a meridian $\mathfrak{m}$ and a longitude $\mathfrak{l}$ for $K$, and orient them arbitrarily.
Given $n\in\mathbb{Z}$, consider the Dehn filling $M_n:=E_K\cup_{{\rm gl}_n}(D^2\times S^1),$
where
$${\rm gl}_n:S^1\times S^1=\partial(D^2\times S^1)\to \partial E_K=S^1\times S^1$$
is a map sending $S^1\times\{1\}$ to a loop homotopic to $\mathfrak{m}^{-n}\mathfrak{l}$.
Then $H_1(M_n;\mathbb{Z})\cong\mathbb{Z}/n\mathbb{Z}$, so $M_n$ is a rational homology sphere.
By Thurston's result (see \cite[Theorem 5.8.2]{Th80}), $M_n$ is hyperbolic for all but finitely many $n$.

Let $e_n$ denote the image of $e$ under the morphism $\mathcal{S}(E_K)\to \mathcal{S}(M_n)$ induced by the inclusion $E_K\subset M_n$.
Then $(q^2-1)e_n=0$ in $\mathcal{S}(M_n)$.

\begin{prop}
Except for finitely many $n$, the element $e_n$ does not vanish in $\mathcal{S}(M_n)$, so $e_n$ is a torsion.
\end{prop}

\begin{proof}
Let $\mathcal{R}^\ast$ denote the set of representations $\sigma:\pi_1(E_K)\to{\rm SL}(2,\mathbb{C})$ such that $\sigma(\tilde{\mathfrak{x}}_1)\sigma(\mathfrak{x}_1)=\mathbf{e}$ and $\sigma(\mathfrak{m}),\sigma(\mathfrak{l})$ are both upper-triangular.
Each character in $\mathcal{X}^\ast$ can be realized by such a representation.
For $\sigma\in\mathcal{R}^\ast$, let $\mu_\sigma$, $\lambda_\sigma$ respectively denote the upper-left entries of $\sigma(\mathfrak{m})$ and $\sigma(\mathfrak{l})$.

Clearly, $\sigma\in\mathcal{R}^\ast$ induces a representation $\pi_1(M_n)\to{\rm SL}(2,\mathbb{C})$ if and only if $\lambda_\sigma=\mu_\sigma^{n}$. So $\mathcal{X}(M_n)$ contains a subset which identifies with
$$\mathcal{Z}_n:=\{\chi_\sigma\colon \sigma\in\mathcal{R}^\ast, \lambda_\sigma=\mu_\sigma^n\}.$$

Let $\mathcal{V}$ denote the Zariski closure of $\mathcal{V}_0:=\{(\lambda_\sigma,\mu_\sigma)\colon\sigma\in\mathcal{R}^\ast\}\subset\mathbb{C}^2$.
It is known that each component of $\mathcal{V}$ has dimension $\le 1$ (see \cite[Page 303]{LR03}, and also \cite[Section 2]{CCGLS94}).
Since $\mu_\sigma$ can take any generic value, there is at least one 1-dimensional component of $\mathcal{V}$ on which $\mu$ is nonconstant.
Let $\mathcal{C}_1,\ldots,\mathcal{C}_r$ be such components. Each $\mathcal{C}_i$ is an affine algebraic curve, defined by some polynomial $A^\ast_i(\lambda,\mu)$ such that the degree in $\lambda$ is positive and $\lambda\nmid A^\ast_i(\lambda,\mu)$.
Note that $\mathcal{C}_i\setminus\mathcal{V}_0$ is a finite set.
Put
$$A^\ast(\lambda,\mu)={\prod}_{i=1}^rA^\ast_i(\lambda,\mu)\in\mathbb{C}[\lambda,\mu].$$
Then $A^\ast$ is a factor of the {\it A-polynomial} of $K$. 
We may write
$$A^\ast(\lambda,\mu)={\sum}_{k=0}^df_k(\mu)\lambda^k, \qquad f_k\in\mathbb{C}[\mu],$$
such that $d\ge 1$, $f_0,f_d\ne 0$.
Remember that $\mu-z\nmid A^\ast(\lambda,\mu)$ for any $z\in\mathbb{C}$.

\begin{cla}\label{claim}
Given any finite set $0\notin Z\subset\mathbb{C}$, when $n$ is sufficiently large, the polynomial $F_n(\mu):={\sum}_{k=0}^df_k(\mu)\mu^{nk}$ has a root in $\mathbb{C}\setminus\{0\}\setminus Z$.
\end{cla}

Consequently, for any prescribed finite set $Y$, when $n$ is sufficiently large,
there exists $\kappa\in\mathbb{C}\setminus \{0,\pm1\}\setminus Y$ such that $F_n(\kappa)=0$, i.e. $A^\ast(\kappa^n,\kappa)=0$.
In particular, we can ensure $(\kappa^n,\kappa)\in\mathcal{V}_0$ and that $\mathsf{t}=\kappa+\kappa^{-1}\ne\pm2$ is generic.

This amounts to saying that there exists $\sigma\in\mathcal{R}^\ast$ with $\mu_\sigma=\kappa$ and $\lambda_\sigma=\kappa^n$.
Fix $\mathsf{t}_{41},\mathsf{t}_{12},\mathsf{t}_{23},\mathsf{t}_{34}$ as in the beginning of this subsection.
Recall that $\mathcal{X}_1$ is 1-dimensional, parameterized by $\mathsf{t}_{24}$. Since $\lambda$ is locally constant when $\mu$ is fixed, we can perturb $\sigma$ into some $\tau\in\mathcal{R}^\ast$ such that $\mu_\tau=\kappa$, $\lambda_\tau=\kappa^n$, and $\mathsf{t}_{24}(\chi_\tau)$ does not satisfy (\ref{eq:quadratic-1}) or (\ref{eq:quadratic-2}). Then $\chi_\tau\in\mathcal{Z}_n$ and $\epsilon(e_n)(\chi_\tau)\ne0$.

In conclusion, except for finitely many $n$, one has $\epsilon(e_n)\ne0$ as a function on $\mathcal{Z}_n$, so that $\epsilon(e_n)\ne0$ on $\mathcal{X}(M_n)$, implying $e_n\ne 0$ in $\mathcal{S}(M_n)$.

\end{proof}

\begin{proof}[Proof of Claim \ref{claim}]
Assume all the roots of $F_n$ are in $\{0\}\cup Z$. Then there exist $a,b_z\in\mathbb{Z}_{\ge 0}$, $C\ne 0$ such that
\begin{align}
F_n(\mu)={\sum}_{k=0}^df_k(\mu)\mu^{nk}=C\mu^a{\prod}_{z\in Z}(\mu-z)^{b_z}.  \label{eq:assumption}
\end{align}
If $n$ is sufficiently large, then $a\le\deg f_0$, and $a+\sum_{z\in Z}b_z=nd+\deg f_d$, so there exists $z_0\in Z$ with $b_{z_0}>d$.
Let
$$G_n(\nu)=F_n(z_0\nu)={\sum}_{k=0}^dg_k(\nu)\nu^{nk}, \qquad g_k(\nu)=f_k(z_0\nu)z_0^{nk}.$$
As an immediate consequence of (\ref{eq:assumption}), $G_n(1)=0$ and the $i$-th derivative $G_n^{(i)}(1)=0$ for $1\le i\le d$.
By the Leibniz product rule,
$$G_n^{(i)}(1)=g_0^{(i)}(1)+\sum_{k=1}^d\sum_{s=0}^i{i\choose s}\cdot g_k^{(i-s)}(1)\cdot nk(nk-1)\cdots(nk-s+1).$$
It can be written as a polynomial in $n$ of degree $\le i$; the coefficient of $n^i$ is $\sum_{k=1}^dk^ig_k(1)$, and the other coefficients are determined by the $g_k^{(j)}(1)$'s.

If $\sum_{k=1}^dk^ig_k(1)=0$ for all $1\le i\le d$, then $g_1(1)=\cdots=g_d(1)=0$, which together with $G_n(1)=0$ would imply $g_0(1)=0$ so that $\mu-z_0\mid A^\ast(\lambda,\mu)$.
Hence $\sum_{k=1}^dk^ig_k(1)\ne 0$ for some $i$, implying that $n$ is a root of $G_n^{(i)}(1)$ which is a polynomial of degree $i$ specified by $A^\ast$ and $z_0$.

Thus, (\ref{eq:assumption}) occurs for only finitely many $n$.
\end{proof}

\newpage

\section*{Appendix}

\bigskip

\noindent
Haimiao Chen (orcid: 0000-0001-8194-1264)\ \ \ \ {\it chenhm@math.pku.edu.cn} \\
Department of Mathematics, Beijing Technology and Business University, \\
Liangxiang Higher Education Park, Fangshan District, Beijing, China.

\end{document}